\documentclass[12pt, twoside]{article}

\usepackage{amsmath}
\usepackage{amsfonts}
\usepackage{amssymb}

\newtheorem{thm}{Theorem}[section]
\newtheorem{cor}[thm]{Corollary}
\newtheorem{lem}[thm]{Lemma}
\newtheorem{prop}[thm]{Proposition}

\begin{document}

\title{Observation estimate for kinetic transport \\
equation by diffusion approximation}
\author{Claude Bardos\thanks{
Universit\'{e} Denis Diderot, Laboratoire J.-L. Lions, 4 Place Jussieu,
BP187, 75252 Paris Cedex 05, France. E-mail address: claude.bardos@gmail.com}
, Kim Dang Phung\thanks{
Universit\'{e} d'Orl\'{e}ans, Laboratoire MAPMO, CNRS UMR 7349, F\'{e}d\'{e}%
ration Denis Poisson, FR CNRS 2964, B\^{a}timent de Math\'{e}matiques, B.P.
6759, 45067 Orl\'{e}ans Cedex 2, France. E-mail address:
kim\_dang\_phung@yahoo.fr. }}
\date{}
\maketitle

\begin{abstract}
We study the unique continuation property for the neutron transport equation
and for a simplified model of the Fokker-Planck equation in a bounded domain
with absorbing boundary condition. An observation estimate is derived. It
depends on the smallness of the mean free path and the frequency of the
velocity average of the initial data. The proof relies on the well known
diffusion approximation under convenience scaling and on basic properties of
this diffusion. Eventually we propose a direct proof for the observation at
one time of parabolic equations. It is based on the analysis of the heat
kernel.
\end{abstract}

\bigskip

\bigskip

\section{Introduction}

\bigskip

This article is devoted to the question of unique continuation for linear
kinetic transport equation with a scattering operator in the diffusive
limit. Let $\Omega$ be a bounded open subset of $\mathbb{R}^{\text{d}}$, d$%
>1 $, with boundary $\partial\Omega$ of class $C^{2}$. Consider in $%
\{(x,v)\in \Omega\times\mathbb{S}^{\text{d}-1}\}\times\mathbb{R}_{t}^{+}$
the transport equation in the $v$ direction with a scattering operator $S$
and absorbing boundary condition 
\begin{equation}
\left\{ 
\begin{array}{ll}
{\partial}_{t}f+\dfrac{1}{\epsilon}v\cdot\nabla f+\dfrac{a}{\epsilon^{2}}%
S\left( f\right) =0 & \quad\text{in}~\Omega\times\mathbb{S}^{\text{d}%
-1}\times\left( 0,+\infty\right) \text{ ,} \\ 
f=0 & \quad\text{on}~\left( \partial\Omega\times\mathbb{S}^{\text{d}%
-1}\right) _{-}\times\left( 0,+\infty\right) \text{ ,} \\ 
f\left( \cdot,\cdot,0\right) =f_{0}\in L^{2}(\Omega\times\mathbb{S}^{\text{d}%
-1})\text{ ,} & 
\end{array}
\right.  \tag{1.1}  \label{1.1}
\end{equation}
where $\epsilon\in\left( 0,1\right] $ is a small parameter and $a\in
L^{\infty}\left( \Omega\right) $ is a scattering opacity satisfying $%
0<c_{min}\leq a\left( x\right) \leq c_{max}<\infty$. Here, $\nabla
=\nabla_{x}$ and $\left( \partial\Omega\times\mathbb{S}^{\text{d}-1}\right)
_{-}=\left\{ \left( x,v\right) \in\partial\Omega\times\mathbb{S}^{\text{d}%
-1};v\cdot\vec{n}_{x}<0\right\} $ where $\vec{n}_{x}$ is the unit outward
normal field at $x\in\partial\Omega$.

Two standard examples of scattering operators $S$ are the following:

\begin{itemize}
\item The neutron scattering operator:%
\begin{equation*}
S=f-\langle f\rangle\text{ where }\langle f\rangle\left( x,t\right) =\dfrac{1%
}{\left\vert \mathbb{S}^{\text{d}-1}\right\vert }\displaystyle\int _{\mathbb{%
S}^{\text{d}-1}}f\left( x,v,t\right) dv\text{ .}
\end{equation*}

\item The Fokker-Planck scattering operator:%
\begin{equation*}
S=-\frac{1}{\text{d}-1}\Delta_{\mathbb{S}^{\text{d}-1}}f\text{ where }%
\Delta_{\mathbb{S}^{\text{d}-1}}\text{ is the Laplace-Beltrami operator on }%
\mathbb{S}^{\text{d}-1}\text{.}
\end{equation*}
\end{itemize}

Let $\omega $ be a nonempty open subset of $\Omega $. Suppose we observe the
solution $f$ at time $T>0$ and on $\omega $, i.e. $f\left( x,v,T\right)
_{\left\vert \left( x,v\right) \in \omega \times \mathbb{S}^{\text{d}%
-1}\right. }$ is known. A classical inverse problem consists to recover at
least one solution, and in particular its initial data, which fits the
observation on $\omega \times \mathbb{S}^{\text{d}-1}\times \left\{
T\right\} $. Our problem of unique continuation is: With how many initial
data, the corresponding solution achieves the given observation $f\left(
x,v,T\right) _{\left\vert \left( x,v\right) \in \omega \times \mathbb{S}^{%
\text{d}-1}\right. }$. Here $\epsilon $ is a small parameter and it is
natural to focus on the limit solution. This is the diffusion approximation
saying that the solution $f$ converges to a solution of a parabolic equation
when $\epsilon $ tends to $0$ (see \cite{B},\cite{DL},\cite{LK},\cite{BR},%
\cite{BGPS},\cite{BSS},\cite{BBGS}). In this framework, two remarks are in
order:

\begin{itemize}
\item For our scattering operator, there holds 
\begin{equation*}
\left\Vert f-\langle f\rangle \right\Vert _{L^{2}(\Omega \times \mathbb{S}^{%
\text{d}-1}\times \mathbb{R}_{t}^{+})}\leq \epsilon \frac{1}{\sqrt{2c_{min}}}%
\Vert f_{0}\Vert _{L^{2}(\Omega \times \mathbb{S}^{\text{d}-1})}\text{ .}
\end{equation*}%
For the operator of neutron transport, one uses a standard energy method by
multiplying both sides of the first line of (\ref{1.1}) by $f$ and
integrating over $\Omega \times \mathbb{S}^{\text{d}-1}\times \left(
0,T\right) $. For the Fokker-Planck scattering operator, one combines the
standard energy method as above and Poincar\'{e} inequality%
\begin{equation*}
\begin{array}{ll}
& \quad \left\Vert f-\langle f\rangle \right\Vert _{L^{2}(\Omega \times 
\mathbb{S}^{\text{d}-1}\times \mathbb{R}_{t}^{+})} \\ 
& \leq \frac{1}{\sqrt{\text{d}-1}}\left\Vert \nabla _{\mathbb{S}^{\text{d}%
-1}}f\right\Vert _{L^{2}(\Omega \times \mathbb{S}^{\text{d}-1}\times \mathbb{%
R}_{t}^{+})}\leq \epsilon \frac{1}{\sqrt{2c_{min}}}\Vert f_{0}\Vert
_{L^{2}(\Omega \times \mathbb{S}^{\text{d}-1})}\text{ .}%
\end{array}%
\end{equation*}

\item In the sense of distributions in $\Omega $, for any $t\geq 0$, the
average of $f$ solves the following parabolic equation%
\begin{equation}
\begin{array}{ll}
& \quad {\partial }_{t}\langle f\rangle -\frac{1}{\text{d}}\nabla \cdot
\left( \frac{1}{a}\nabla \langle f\rangle \right) \\ 
& =\nabla \cdot \left( \frac{1}{a}\langle \left( v\otimes v\right) \nabla
\left( f-\langle f\rangle \right) \rangle \right) +\epsilon \nabla \cdot
\left( \frac{1}{a}\langle v\partial _{t}f\rangle \right) \text{ .}%
\end{array}
\tag{1.2}  \label{1.2}
\end{equation}%
Indeed, multiply by $\dfrac{\epsilon }{a}v$ the equation ${\partial }_{t}f+%
\dfrac{1}{\epsilon }v\cdot \nabla f+\dfrac{a}{\epsilon ^{2}}Sf=0$ and take
the average over $\mathbb{S}^{\text{d}-1}$, using ${\partial }_{t}\langle
f\rangle +\dfrac{1}{\epsilon }\langle v\cdot \nabla f\rangle =0$, $\langle
v\cdot \nabla Sf\rangle =\langle v\cdot \nabla f\rangle $ and $\langle
v\left( v\cdot \nabla \langle f\rangle \right) \rangle =\frac{1}{\text{d}}%
\nabla \langle f\rangle $, one obtains for any $t\geq 0$ and any $\varphi
\in C_{0}^{\infty }\left( \Omega \right) $%
\begin{equation*}
\int_{\Omega }{\partial }_{t}\langle f\rangle \varphi dx+\frac{1}{\text{d}}%
\int_{\Omega }\frac{1}{a}\nabla \langle f\rangle \cdot \nabla \varphi
dx+\int_{\Omega }\frac{1}{a}\langle v\left( v\cdot \nabla \left( f-\langle
f\rangle \right) +\epsilon \partial _{t}f\right) \rangle \cdot \nabla
\varphi dx=0\text{ .}
\end{equation*}%
Moreover, we prove that the boundary condition on $\langle f\rangle $ is
small in some adequate norm with respect to $\epsilon $. In the sequel, any
estimates will be explicit with respect to $\epsilon $.
\end{itemize}

\bigskip

Backward uniqueness for parabolic equation has a long history (see \cite{DJP}%
,\cite{V}). Lions and Malgrange \cite{LM} used the method of Carleman
estimates. Later, Bardos and Tartar \cite{BT} gave some improvements by
using the log convexity method of Agmon and Nirenberg. More recently,
motivated by control theory and inverse problems (see \cite{I},\cite{P}),
Carleman estimates became an important tool to achieve an observability
inequality (see \cite{FI},\cite{FZ},\cite{FG},\cite{LRL},\cite{LeRR},\cite%
{LRR}). In \cite{PW}, the desired observability inequality is deduced from
the observation estimate at one point in time which is obtained by studying
the frequency function in the spirit of log convexity method. In particular,
one can quantify the following unique continuation property (see \cite{EFV},%
\cite{PWa}): If $u\left( x,t\right) =e^{t\Delta}u_{0}\left( x\right) $ with $%
u_{0}\in L^{2}\left( \Omega\right) $ and $u\left( \cdot,T\right) =0$ on $%
\omega$, then $u_{0}\equiv0$.

\bigskip

Our main result below involves the regularity of the nonzero initial data $%
f_{0}$ measured in term of two quantities. Let $p>2$,%
\begin{equation*}
\mathbb{M}_{p}:=\frac{\left\Vert f_{0}\right\Vert _{L^{2p}\left( \Omega
\times \mathbb{S}^{\text{d}-1}\right) }}{\left\Vert \langle f_{0}\rangle
\right\Vert _{L^{2}\left( \Omega \right) }}\text{ and }\mathbb{F}:=\frac{%
\left\Vert \langle f_{0}\rangle \right\Vert _{L^{2}\left( \Omega \right)
}^{2}}{\left\Vert \langle f_{0}\rangle \right\Vert _{H^{-1}\left( \Omega
\right) }^{2}}\text{ .}
\end{equation*}%
Observe in particular that $\mathbb{F}$ is the most natural evaluation of
the frequency of the velocity average of the initial data.

\bigskip

\begin{thm}
\label{theorem1.1} Suppose that $a\in C^{2}\left( \overline{\Omega}\right) $
and $f_{0}\in L^{2p}\left( \Omega\times\mathbb{S}^{\normalfont{\text{d}}%
-1}\right) $ with $M_{p}+F<+\infty$ for some $p>2$. Then the unique solution 
$f$ of (\ref{1.1}) satisfies for any $T>0$ 
\begin{equation*}
\left( 1-\epsilon^{\frac{1}{2p}}(1+T^{\frac{p-1}{2p}}C_{p})\mathbb{M}%
_{p}e^{\sigma\left( f_{0},T\right) }\right) \left\Vert \langle
f_{0}\rangle\right\Vert _{L^{2}\left( \Omega\right) }\leq e^{\sigma\left(
f_{0},T\right) }\left\Vert \langle f\rangle\left( \cdot,T\right) \right\Vert
_{L^{2}\left( \omega\right) }
\end{equation*}
with $C_{p}=\left( \frac{p-1}{p-2}\right) ^{\frac{p-1}{2p}}\left( \frac {1}{p%
}\right) ^{\frac{1}{2p}}$ and $\sigma\left( f_{0},T\right) =c\left( 1+\frac{1%
}{T}+T\mathbb{F}\right) $ where $c$ only depends on $\left( \Omega,\omega,%
\normalfont{\text{d}},a\right) $.
\end{thm}

\bigskip

By a direct application of our main result, we have:

\bigskip

\begin{cor}
\label{corollary1.2} Let $a\in C^{2}\left( \overline{\Omega }\right) $ and $%
f_{0}\in L^{2p}\left( \Omega \times \mathbb{S}^{\normalfont{\text{d}}%
-1}\right) $ with $M_{p}+F<+\infty $ for some $p>2$. Suppose that $f_{0}\geq
0$. Then there is $\epsilon _{0}\in \left( 0,1\right) $ depending on $\left( 
\mathbb{M}_{p},\mathbb{F},\Omega ,\omega ,\normalfont{\text{d}},p,T,a\right) 
$ such that if $f\left( \cdot ,\cdot ,T\right) =0$ on $\omega \times \mathbb{%
S}^{\normalfont{\text{d}}-1}$ for some $\epsilon \leq \epsilon _{0}$, then $%
f_{0}\equiv 0$.
\end{cor}

\bigskip

\bigskip

This paper is organized as follows: The proof of the main result is given in
the next section. It requires two important results: an approximation
diffusion convergence of the average of $f$; an observation estimate at one
point in time for the diffusion equation with homogeneous Dirichlet boundary
condition. In Section 3, we prove the approximation theorem stated in
Section 2. In Section 4, a direct proof of the observation inequality at one
point in time for parabolic equations is proposed. Finally, in an appendix,
we prove a backward estimate for the diffusion equation and a trace estimate
for the kinetic transport equation.

\bigskip

\bigskip

\section{Proof of main Theorem \protect\ref{theorem1.1}}

\bigskip

The main task in the proof of Theorem \ref{theorem1.1} consists on the two
following propositions. Below we denote by $u\in C\left( \left[ 0,T\right]
,L^{2}(\Omega )\right) \cap L^{2}\left( 0,T;H_{0}^{1}\left( \Omega \right)
\right) $ any solution of the diffusion equation%
\begin{equation}
{\partial }_{t}u-\frac{1}{\text{d}}\nabla \cdot \left( \frac{1}{a}\nabla
u\right) =0  \tag{2.1}  \label{2.1}
\end{equation}%
with $a\in C^{2}\left( \overline{\Omega }\right) $ and $0<c_{min}\leq
a\left( x\right) \leq c_{max}<\infty $.

\bigskip

\begin{prop}
\label{proposition2.1} There are $C>0$\ and $\mu \in \left( 0,1\right) $\
such that any solution$\ u\in C\left( \left[ 0,T\right] ,L^{2}(\Omega
)\right) \cap L^{2}\left( 0,T;H_{0}^{1}\left( \Omega \right) \right) $ of (%
\ref{2.1}) satisfies 
\begin{equation*}
\int_{\Omega }\left\vert u\left( x,T\right) \right\vert ^{2}dx\leq \left(
Ce^{\frac{C}{T}}\int_{\omega }\left\vert u\left( x,T\right) \right\vert
^{2}dx\right) ^{1-\mu }\left( \int_{\Omega }\left\vert u\left( x,0\right)
\right\vert ^{2}dx\right) ^{\mu }\text{ .}
\end{equation*}%
Here $C$\ and $\mu $\ only depend on $\left( a,\Omega ,\omega ,%
\normalfont{\text{d}}\right) $.
\end{prop}

\bigskip

As an immediate application, combining with the following backward estimate
for diffusion equation%
\begin{equation}
\left\Vert u\left( \cdot ,0\right) \right\Vert _{L^{2}\left( \Omega \right)
}\leq ce^{cT\frac{\left\Vert u\left( \cdot ,0\right) \right\Vert
_{L^{2}\left( \Omega \right) }^{2}}{\left\Vert u\left( \cdot ,0\right)
\right\Vert _{H^{-1}\left( \Omega \right) }^{2}}}\left\Vert u\left( \cdot
,T\right) \right\Vert _{L^{2}\left( \Omega \right) } \text{ ,}  \tag{2.2}
\label{2.2}
\end{equation}%
we have:

\bigskip

\begin{cor}
\label{corollary2.2} For any nonzero $u\in C\left( \left[ 0,T\right]
,L^{2}(\Omega )\right) \cap L^{2}\left( 0,T;H_{0}^{1}\left( \Omega \right)
\right) $ solution of (\ref{2.1}), one has 
\begin{equation*}
\left\Vert u\left( \cdot ,0\right) \right\Vert _{L^{2}\left( \Omega \right)
}\leq e^{C\left( 1+\frac{1}{T}+T\frac{\left\Vert u\left( \cdot ,0\right)
\right\Vert _{L^{2}\left( \Omega \right) }^{2}}{\left\Vert u\left( \cdot
,0\right) \right\Vert _{H^{-1}\left( \Omega \right) }^{2}}\right)
}\left\Vert u\left( \cdot ,T\right) \right\Vert _{L^{2}\left( \omega \right)
}
\end{equation*}%
where $C$\ only depends on $\left( a,\Omega ,\omega ,\normalfont{\text{d}}%
\right) $.
\end{cor}

\bigskip

\bigskip

\begin{prop}
\label{proposition2.3} Assume $f_{0}\in L^{2p}(\Omega \times \mathbb{S}^{%
\normalfont{\text{d}}-1})$ for some $p>2$ and consider $u\in C\left(
0,T;H_{0}^{1}(\Omega )\right) $ solution of (\ref{2.1}) with initial data $%
u\left( \cdot ,0\right) =\langle f_{0}\rangle $, then for any $T>0$ and any $%
\chi \in C_{0}^{\infty }\left( \omega \right) $, the solution $f$ of (\ref%
{1.1}) satisfies 
\begin{equation*}
\left\Vert \chi \left( \langle f\rangle _{\left\vert t=T\right. }-u\left(
\cdot ,T\right) \right) \right\Vert _{H^{-1}\left( \Omega \right) }\leq
\epsilon ^{\frac{1}{2p}}\left( 1+T^{\frac{p-1}{2p}}C_{p}\right) C\Vert
f_{0}\Vert _{L^{2p}(\Omega \times \mathbb{S}^{\normalfont{\text{d}}-1})}
\end{equation*}%
where $C_{p}=\left( \frac{p-1}{p-2}\right) ^{\frac{p-1}{2p}}\left( \frac{1}{p%
}\right) ^{\frac{1}{2p}}$ and $C>0$ only depends on $\left( \Omega ,%
\normalfont{\text{d}},a,\chi \right) $.
\end{prop}

\bigskip

The proof of Proposition \ref{proposition2.1} and Proposition \ref%
{proposition2.3} is given in section 4 and section 3 respectively.

\bigskip

\noindent In one hand, since $u\left( \cdot ,0\right) =\langle f_{0}\rangle $%
, we have by Corollary \ref{corollary2.2}%
\begin{equation*}
\left\Vert \langle f_{0}\rangle \right\Vert _{L^{2}\left( \Omega \right)
}\leq e^{C\left( 1+\frac{1}{T}+T\frac{\left\Vert \langle f_{0}\rangle
\right\Vert _{L^{2}\left( \Omega \right) }^{2}}{\left\Vert \langle
f_{0}\rangle \right\Vert _{H^{-1}\left( \Omega \right) }^{2}}\right)
}\left\Vert \chi u\left( \cdot ,T\right) \right\Vert _{L^{2}\left( \Omega
\right) }\text{ .}
\end{equation*}%
On the other hand, by regularizing effect, we have%
\begin{equation*}
\begin{array}{ll}
\left\Vert \chi u\left( \cdot ,T\right) \right\Vert _{L^{2}\left( \Omega
\right) } & \leq \left\Vert \chi u\left( \cdot ,T\right) \right\Vert
_{H^{-1}\left( \Omega \right) }^{1/2}\left\Vert \chi u\left( \cdot ,T\right)
\right\Vert _{H_{0}^{1}\left( \Omega \right) }^{1/2} \\ 
& \leq C\left\Vert \chi u\left( \cdot ,T\right) \right\Vert _{H^{-1}\left(
\Omega \right) }^{1/2}\left( 1+\frac{1}{T^{1/4}}\right) \left\Vert \langle
f_{0}\rangle \right\Vert _{L^{2}\left( \Omega \right) }^{1/2}\text{ .}%
\end{array}%
\end{equation*}%
Therefore, the two above facts yield 
\begin{equation*}
\begin{array}{ll}
\left\Vert \langle f_{0}\rangle \right\Vert _{L^{2}\left( \Omega \right) } & 
\leq e^{C\left( 1+\frac{1}{T}+T\frac{\left\Vert \langle f_{0}\rangle
\right\Vert _{L^{2}\left( \Omega \right) }^{2}}{\left\Vert \langle
f_{0}\rangle \right\Vert _{H^{-1}\left( \Omega \right) }^{2}}\right) }\left(
\left\Vert \chi \left( u\left( \cdot ,T\right) -\langle f\rangle
_{\left\vert t=T\right. }\right) \right\Vert _{H^{-1}\left( \Omega \right)
}+\left\Vert \chi \langle f\rangle _{\left\vert t=T\right. }\right\Vert
_{H^{-1}\left( \Omega \right) }\right) \\ 
& \leq e^{C\left( 1+\frac{1}{T}+T\frac{\left\Vert \langle f_{0}\rangle
\right\Vert _{L^{2}\left( \Omega \right) }^{2}}{\left\Vert \langle
f_{0}\rangle \right\Vert _{H^{-1}\left( \Omega \right) }^{2}}\right) }\left(
\epsilon ^{\frac{1}{2p}}\left( 1+T^{\frac{p-1}{2p}}C_{p}\right) \Vert
f_{0}\Vert _{L^{2p}(\Omega \times \mathbb{S}^{\text{d}-1})}+\left\Vert \chi
\langle f\rangle _{\left\vert t=T\right. }\right\Vert _{H^{-1}\left( \Omega
\right) }\right)%
\end{array}%
\end{equation*}%
where in the last line we used Proposition \ref{proposition2.3}. This
completes the proof.

\bigskip

\bigskip

\section{Estimates for diffusion approximation}

\bigskip

Below we give precise error estimates for the diffusion approximation.

\bigskip

\begin{thm}
\label{theorem3.1} Let $a\in C^{1}\left( \overline{\Omega }\right) $ such
that $0<c_{min}\leq a\left( x\right) \leq c_{max}<\infty $. Assume $f_{0}\in
L^{2p}(\Omega \times \mathbb{S}^{\normalfont{\text{d}}-1})$ for some $p>2$
and consider $u\in C\left( 0,T;H_{0}^{1}(\Omega )\right) $ solution of (\ref%
{2.1}) with initial data $u\left( \cdot ,0\right) =\langle f_{0}\rangle $,
then for any $T>0$, the solution $f$ of (\ref{1.1}) satisfies%
\begin{equation*}
\left\Vert \langle f\rangle _{\left\vert t=T\right. }-u\left( \cdot
,T\right) \right\Vert _{H^{-1}\left( \Omega \right) }+\left\Vert \langle
f\rangle -u\right\Vert _{L^{2}\left( \Omega \times \left( 0,T\right) \right)
}\leq \epsilon ^{\frac{1}{2p}}\left( 1+T^{\frac{p-1}{2p}}C_{p}\right) C\Vert
f_{0}\Vert _{L^{2p}(\Omega \times \mathbb{S}^{\normalfont{\text{d}}-1})}
\end{equation*}%
where $C_{p}=\left( \frac{p-1}{p-2}\right) ^{\frac{p-1}{2p}}\left( \frac{1}{p%
}\right) ^{\frac{1}{2p}}$ and $C>0$ only depends on $\left( \Omega ,%
\normalfont{\text{d}}\right) $ and $\left( c_{min},c_{max},\left\Vert \nabla
a\right\Vert _{\infty }\right) $.
\end{thm}

\bigskip

\bigskip

In the literature, there are at least two ways to get diffusion
approximation estimates:

- Use a Hilbert expansion: The solution $f$ of the transport problem can be
formally written as $f=f_{0}+\epsilon f_{1}+\epsilon^{2}f_{2}+...$ and we
substitute this expansion into the governing equations in order to prove
existence of $f_{0},f_{1},f_{2},...$. Next we set $F=f-\left( f_{0}+\epsilon
f_{1}\right) $ and check that it solves a transport problem for which energy
method can be used. This way requires well-prepared initial data that is $%
f_{0}=\langle f_{0}\rangle$ to avoid initial layers.

- Use moment method: The zeroth and first moments of $f$ are respectively $%
\langle f\rangle $ and $\langle vf\rangle $. First, we check that $f-\langle
f\rangle $ is small in some adequate norm with respect to $\epsilon $. Next
by computing the zeroth and first moments of the equation solved by $f$ (as
it was done in the introduction), we derive that $\langle f\rangle $ solves
a parabolic problem for which energy method can be used. This way and a new $%
\epsilon $ uniform estimate on the trace (see Proposition \ref%
{proposition3.2} below) give Theorem \ref{theorem3.1}. Notice that since
only the average of $f$, is involved, the proof requires no analysis of the
initial layer near $t=0$.

\bigskip

\begin{prop}
\label{proposition3.2} If $f_{0}\in L^{2p}(\Omega \times \mathbb{S}^{%
\normalfont{\text{d}}-1})$ for some $p>2$, then the solution $f$ of (\ref%
{1.1}) satisfies 
\begin{equation*}
\left\Vert f\right\Vert _{L^{2}\left( \partial \Omega \times \mathbb{S}^{%
\normalfont{\text{d}}-1}\times \left( 0,T\right) \right) }\leq CT^{\frac{p-1%
}{2p}}\epsilon ^{\frac{1}{2p}}C_{p}\Vert f_{0}\Vert _{L^{2p}(\Omega \times 
\mathbb{S}^{\text{d}-1})}
\end{equation*}%
where $C_{p}=\left( \frac{p-1}{p-2}\right) ^{\frac{p-1}{2p}}\left( \frac{1}{p%
}\right) ^{\frac{1}{2p}}$ and $C>0$ only depends on $\left( \Omega ,%
\normalfont{\text{d}}\right) $.
\end{prop}

\bigskip

Proposition \ref{proposition3.2} is proved in Appendix. The proof of Theorem %
\ref{theorem3.1} starts as follows. Let $w_{\epsilon }=\langle f\rangle -u$
where $u$ solves 
\begin{equation*}
\left\{ 
\begin{array}{ll}
{\partial }_{t}u-\frac{1}{\text{d}}\nabla \cdot \left( \frac{1}{a}\nabla
u\right) =0 & \quad \text{in}~\Omega \times \left( 0,+\infty \right) \text{ ,%
} \\ 
u=0 & \quad \text{on}~\partial \Omega \times \left( 0,+\infty \right) \text{
,} \\ 
u\left( \cdot ,0\right) =\langle f_{0}\rangle \in L^{2}(\Omega )\text{ .} & 
\end{array}%
\right.
\end{equation*}%
By (\ref{1.2}) and a density argument, $w_{\epsilon }$ solves for any $t\geq
0$ and any $\varphi \in H_{0}^{1}\left( \Omega \right) $%
\begin{equation}
\begin{array}{ll}
& \quad \displaystyle\int_{\Omega }{\partial }_{t}w_{\epsilon }\varphi dx+%
\displaystyle\frac{1}{\text{d}}\int_{\Omega }\nabla w_{\epsilon }\cdot \frac{%
1}{a}\nabla \varphi dx \\ 
& =-\displaystyle\int_{\Omega }\langle v\left( v\cdot \nabla \left(
f-\langle f\rangle \right) \right) \rangle \cdot \frac{1}{a}\nabla \varphi
dx-\displaystyle\epsilon \int_{\Omega }\langle v\partial _{t}f\rangle \cdot 
\frac{1}{a}\nabla \varphi dx%
\end{array}
\tag{3.1}  \label{3.1}
\end{equation}%
with boundary condition $w_{\epsilon }=\langle f\rangle $ on $\partial
\Omega \times \mathbb{R}_{t}^{+}$ and initial data $w_{\epsilon }\left(
\cdot ,0\right) =0$. We choose 
\begin{equation*}
\varphi =\left( -\frac{1}{\text{d}}\nabla \cdot \left( \frac{1}{a}\nabla
\right) \right) ^{-1}w_{\epsilon }\text{ .}
\end{equation*}%
By integrations by parts, the identity (\ref{3.1}) becomes:%
\begin{equation}
\begin{array}{ll}
& \quad \displaystyle\dfrac{1}{2\text{d}}\frac{d}{dt}\int_{\Omega }\frac{1}{a%
}\left\vert \nabla \varphi \right\vert ^{2}dx+\left\Vert \frac{1}{\text{d}}%
\nabla \cdot \left( \frac{1}{a}\nabla \right) \varphi \right\Vert
_{L^{2}\left( \Omega \right) }^{2} \\ 
& =-\displaystyle\int_{\partial \Omega }\langle f\rangle \frac{1}{a}\partial
_{n}\varphi dx \\ 
& \quad -\displaystyle\int_{\Omega }\langle v\left( v\cdot \nabla \left(
f-\langle f\rangle \right) \right) \rangle \cdot \frac{1}{a}\nabla \varphi dx
\\ 
& \quad -\displaystyle\epsilon \int_{\Omega }\langle v\partial _{t}f\rangle
\cdot \frac{1}{a}\nabla \varphi dx\text{ .}%
\end{array}
\tag{3.2}  \label{3.2}
\end{equation}%
First, the contribution of the boundary data is estimate: One has, by a
classical trace theorem 
\begin{equation*}
-\int_{\partial \Omega }\langle f\rangle \frac{1}{a}\partial _{n}\varphi
dx\leq C_{1}\left\Vert f\right\Vert _{L^{2}\left( \partial \Omega \times 
\mathbb{S}^{\text{d}-1}\right) }\left\Vert \nabla \cdot \left( \frac{1}{a}%
\nabla \right) \varphi \right\Vert _{L^{2}\left( \Omega \right) }
\end{equation*}%
where the constant $C_{1}$ depends on $\left\Vert \nabla a\right\Vert
_{\infty }$.

Secondly, the contribution of the term 
\begin{equation*}
\int_{\Omega}\langle v\left( v\cdot\nabla\left( f-\langle f\rangle\right)
\right) \rangle\cdot\frac{1}{a}\nabla\varphi dx
\end{equation*}
is estimated: By integration by parts and using $\nabla\varphi=\partial
_{n}\varphi\vec{n}_{x}$ on $\partial\Omega$, one has 
\begin{equation*}
\begin{array}{ll}
\displaystyle\int_{\Omega}\langle v\left( v\cdot\nabla\left( f-\langle
f\rangle\right) \right) \rangle\cdot\frac{1}{a}\nabla\varphi dx & =-%
\displaystyle\dfrac{1}{\left\vert \mathbb{S}^{\text{d}-1}\right\vert }%
\int_{\Omega\times\mathbb{S}^{\text{d}-1}}(f-\langle f\rangle)v\cdot
\nabla\left( v\cdot\frac{1}{a}\nabla\varphi\right) dxdv \\ 
& \quad+\displaystyle\dfrac{1}{\left\vert \mathbb{S}^{\text{d}-1}\right\vert 
}\int_{\partial\Omega\times\mathbb{S}^{\text{d}-1}}\left( v\cdot\vec{n}%
_{x}\right) ^{2}(f-\langle f\rangle)\frac{1}{a}\partial_{n}\varphi dxdv%
\end{array}%
\end{equation*}
which implies%
\begin{equation*}
\begin{array}{ll}
\displaystyle\int_{\Omega}\langle v\left( v\cdot\nabla\left( f-\langle
f\rangle\right) \right) \rangle\cdot\frac{1}{a}\nabla\varphi dx & \leq
C_{1}\left\Vert f-\langle f\rangle\right\Vert _{L^{2}\left( \Omega \times%
\mathbb{S}^{\text{d}-1}\right) }\left\Vert \nabla\cdot\left( \frac {1}{a}%
\nabla\right) \varphi\right\Vert _{L^{2}\left( \Omega\right) } \\ 
& \quad+C_{1}\left\Vert f\right\Vert _{L^{2}\left( \partial\Omega \times%
\mathbb{S}^{\text{d}-1}\right) }\left\Vert \nabla\cdot\left( \frac {1}{a}%
\nabla\right) \varphi\right\Vert _{L^{2}\left( \Omega\right) }%
\end{array}%
\end{equation*}
with some constant $C_{1}>0$ depending on $\left\Vert \nabla a\right\Vert
_{\infty}$.

Thirdly, the contribution of the term $\displaystyle\epsilon \int_{\Omega
}\langle v\partial _{t}f\rangle \cdot \frac{1}{a}\nabla \varphi dx$ is
estimated: From the identities%
\begin{equation*}
\begin{array}{ll}
& \quad \displaystyle\epsilon \int_{\Omega }\langle v\partial _{t}f\rangle
\cdot \frac{1}{a}\nabla \varphi dx \\ 
& =\displaystyle\frac{1}{\left\vert \mathbb{S}^{\text{d}-1}\right\vert }%
\epsilon \frac{d}{dt}\int_{\Omega \times \mathbb{S}^{\text{d}-1}}fv\cdot 
\frac{1}{a}\nabla \varphi dxdv-\displaystyle\frac{1}{\left\vert \mathbb{S}^{%
\text{d}-1}\right\vert }\int_{\Omega \times \mathbb{S}^{\text{d}-1}}fv\cdot 
\frac{1}{a}\nabla \left( \epsilon {\partial }_{t}\varphi \right) dxdvdt \\ 
& =\displaystyle\frac{1}{\left\vert \mathbb{S}^{\text{d}-1}\right\vert }%
\epsilon \frac{d}{dt}\int_{\Omega \times \mathbb{S}^{\text{d}-1}}fv\cdot 
\frac{1}{a}\nabla \varphi dxdv-\displaystyle\frac{1}{\left\vert \mathbb{S}^{%
\text{d}-1}\right\vert }\int_{\Omega \times \mathbb{S}^{\text{d}-1}}\left(
f-\langle f\rangle \right) v\cdot \frac{1}{a}\nabla \left( \epsilon {%
\partial }_{t}\varphi \right) dxdv%
\end{array}%
\end{equation*}%
and 
\begin{equation*}
\begin{array}{ll}
\epsilon {\partial }_{t}\varphi  & =\left( -\frac{1}{\text{d}}\nabla \cdot
\left( \frac{1}{a}\nabla \right) \right) ^{-1}\left( \epsilon {\partial }%
_{t}w_{\epsilon }\right) =\left( -\frac{1}{\text{d}}\nabla \cdot \left( 
\frac{1}{a}\nabla \right) \right) ^{-1}\left( -\left\langle v\cdot \nabla
f\right\rangle -\epsilon {\partial }_{t}u\right)  \\ 
& =\left( -\frac{1}{\text{d}}\nabla \cdot \left( \frac{1}{a}\nabla \right)
\right) ^{-1}\left\langle -v\cdot \nabla \left( f-\langle f\rangle \right)
\right\rangle +\epsilon u\text{ ,}%
\end{array}%
\end{equation*}%
we see that%
\begin{equation*}
\begin{array}{ll}
& \displaystyle\epsilon \int_{\Omega }\langle v\partial _{t}f\rangle \cdot 
\frac{1}{a}\nabla \varphi dx \\ 
& =\displaystyle\frac{1}{\left\vert \mathbb{S}^{\text{d}-1}\right\vert }%
\epsilon \frac{d}{dt}\int_{\Omega \times \mathbb{S}^{\text{d}-1}}fv\cdot 
\frac{1}{a}\nabla \varphi dxdv \\ 
& \quad +\displaystyle\frac{1}{\left\vert \mathbb{S}^{\text{d}-1}\right\vert 
}\int_{\Omega \times \mathbb{S}^{\text{d}-1}}\left( f-\langle f\rangle
\right) v\cdot \frac{1}{a}\nabla \left( \left( -\frac{1}{\text{d}}\nabla
\cdot \left( \frac{1}{a}\nabla \right) \right) ^{-1}\left\langle v\cdot
\nabla \left( f-\langle f\rangle \right) \right\rangle \right) dxdv \\ 
& \quad -\epsilon \displaystyle\frac{1}{\left\vert \mathbb{S}^{\text{d}%
-1}\right\vert }\int_{\Omega \times \mathbb{S}^{\text{d}-1}}\left( f-\langle
f\rangle \right) v\cdot \frac{1}{a}\nabla udxdv \\ 
& \leq \displaystyle\frac{1}{\left\vert \mathbb{S}^{\text{d}-1}\right\vert }%
\epsilon \frac{d}{dt}\int_{\Omega \times \mathbb{S}^{\text{d}-1}}fv\cdot 
\frac{1}{a}\nabla \varphi dxdv+C\left\Vert f-\langle f\rangle \right\Vert
_{L^{2}\left( \Omega \times \mathbb{S}^{\text{d}-1}\right) }^{2}+\epsilon
^{2}C\left\Vert \nabla u\right\Vert _{L^{2}\left( \Omega \right) }^{2}\text{
.}%
\end{array}%
\end{equation*}%
Combining the three above contributions with (\ref{3.2}), one obtains%
\begin{equation*}
\begin{array}{ll}
\displaystyle\frac{d}{dt}\int_{\Omega }\frac{1}{a}\left\vert \nabla \varphi
\right\vert ^{2}dx+\left\Vert \nabla \cdot \left( \frac{1}{a}\nabla \right)
\varphi \right\Vert _{L^{2}\left( \Omega \right) }^{2} & \leq \displaystyle%
\epsilon C\frac{d}{dt}\int_{\Omega \times \mathbb{S}^{\text{d}-1}}fv\cdot 
\frac{1}{a}\nabla \varphi dxdv+\epsilon ^{2}C\left\Vert \nabla u\right\Vert
_{L^{2}\left( \Omega \right) }^{2} \\ 
& \quad +C\left( \left\Vert f\right\Vert _{L^{2}\left( \partial \Omega
\times \mathbb{S}^{\text{d}-1}\right) }^{2}+\left\Vert f-\langle f\rangle
\right\Vert _{L^{2}\left( \Omega \times \mathbb{S}^{\text{d}-1}\right)
}^{2}\right) \text{ .}%
\end{array}%
\end{equation*}%
Integrating the above over $\left( 0,T\right) $, we observe with $\varphi
=\left( -\frac{1}{\text{d}}\nabla \cdot \left( \frac{1}{a}\nabla \right)
\right) ^{-1}w_{\epsilon }$ and $w_{\epsilon }=\langle f\rangle -u$ that 
\begin{equation*}
\begin{array}{ll}
& \quad \left\Vert w_{\epsilon }(\cdot ,T)\right\Vert _{H^{-1}\left( \Omega
\right) }^{2}+\left\Vert w_{\epsilon }\right\Vert _{L^{2}\left( \Omega
\times \left( 0,T\right) \right) }^{2} \\ 
& \leq \epsilon C\left( \left\Vert f_{\left\vert t=T\right. }\right\Vert
_{L^{2}\left( \Omega \times \mathbb{S}^{\text{d}-1}\right) }^{2}+\Vert
u_{\left\vert t=T\right. }\Vert _{L^{2}(\Omega )}^{2}+\left\Vert
f_{0}\right\Vert _{L^{2}\left( \Omega \times \mathbb{S}^{\text{d}-1}\right)
}^{2}+\Vert u_{0}\Vert _{L^{2}(\Omega )}^{2}\right)  \\ 
& \quad +\epsilon ^{2}C\left\Vert \nabla u\right\Vert _{L^{2}\left( \Omega
\times \left( 0,T\right) \right) }^{2} \\ 
& \quad +C\left( \left\Vert f\right\Vert _{L^{2}\left( \partial \Omega
\times \mathbb{S}^{\text{d}-1}\times \left( 0,T\right) \right)
}^{2}+\left\Vert f-\langle f\rangle \right\Vert _{L^{2}\left( \Omega \times 
\mathbb{S}^{\text{d}-1}\times \left( 0,T\right) \right) }^{2}\right) \text{ .%
}%
\end{array}%
\end{equation*}%
Next, we use the trace estimate in Proposition \ref{proposition3.2}, 
\begin{equation*}
\int_{\Omega \times \mathbb{S}^{\text{d}-1}}\left\vert f\left( x,v,T\right)
\right\vert ^{2}dxdv+\frac{2c_{min}}{\epsilon ^{2}}\int_{0}^{T}\int_{\Omega
\times \mathbb{S}^{\text{d}-1}}\left\vert f-\langle f\rangle \right\vert
^{2}dxdvdt\leq \int_{\Omega \times \mathbb{S}^{\text{d}-1}}\left\vert
f_{0}\right\vert ^{2}dxdv\text{ }
\end{equation*}%
and 
\begin{equation*}
\int_{\Omega }\left\vert u\left( x,T\right) \right\vert ^{2}dx+\frac{2}{%
\text{d}c_{max}}\int_{0}^{T}\int_{\Omega }\left\vert \nabla u\right\vert
^{2}dxdt\leq \int_{\Omega }\left\vert \langle f_{0}\rangle \right\vert ^{2}dx%
\text{ }
\end{equation*}%
to get that 
\begin{equation*}
\begin{array}{ll}
& \quad \left\Vert \left( \langle f\rangle -u\right) \left( \cdot ,T\right)
\right\Vert _{H^{-1}\left( \Omega \right) }+\left\Vert \langle f\rangle
-u\right\Vert _{L^{2}\left( \Omega \times \left( 0,T\right) \right) } \\ 
& \leq \sqrt{\epsilon }C\left\Vert f_{0}\right\Vert _{L^{2}\left( \Omega
\times \mathbb{S}^{\text{d}-1}\right) }+\epsilon ^{\frac{1}{2p}}T^{\frac{p-1%
}{2p}}C_{p}C\Vert f_{0}\Vert _{L^{2p}(\Omega \times \mathbb{S}^{\text{d}-1})}%
\text{ .}%
\end{array}%
\end{equation*}%
This completes the proof.

\bigskip

\bigskip

\section{Observation estimates for diffusion equation}

\bigskip

In this section, we establish an observation estimate at one point in time
for parabolic equations (see Theorem \ref{theorem4.1} below). Such estimate
is an interpolation inequality. H\"{o}lder type inequalities of such form
already appear in \cite{LR} for elliptic operators by Carleman inequalities.
It applies to the observability for the heat equation in manifold and to the
sum of eigenfunctions estimate of Lebeau-Robbiano. On the other hand, for
parabolic operators, Escauriaza, Fernandez and Vessella proved such
interpolation estimate far from the boundary by some adequate Carleman
estimates \cite{EFV}. Here our approach is completely new and uses
properties of the heat kernel with a parametrix of order $0$.

\bigskip

\bigskip

\begin{thm}
\label{theorem4.1} Let $\Omega $\ be a bounded open set in $\mathbb{R}^{n}$, 
$n\geq 1$, either convex or $C^{2}$\ and connected. Let $\omega $\ be a
nonempty open subset of $\Omega $, and $T>0$. Let $A$ be a $n\times n$
symmetric positive-definite matrix with $C^{1}\left( \overline{\Omega }%
\times \left[ 0,T\right] \right) $ coefficients such that $A\left( \cdot
,T\right) \in C^{2}\left( \overline{\Omega }\right) $. There are $c>0$\ and $%
\mu \in \left( 0,1\right) $\ such that any solution to 
\begin{equation*}
\left\{ 
\begin{array}{ll}
\partial _{t}u-\nabla \cdot \left( A\nabla u\right) =0 & \text{\textit{in}~}%
\Omega \times \left( 0,T\right) \ \text{,} \\ 
u=0 & \text{\textit{on}~}\partial \Omega \times \left( 0,T\right) \text{ ,}
\\ 
u\left( \cdot ,0\right) \in L^{2}\left( \Omega \right) \text{ ,} & 
\end{array}%
\right.
\end{equation*}%
satisfies 
\begin{equation*}
\int_{\Omega }\left\vert u\left( x,T\right) \right\vert ^{2}dx\leq \left(
c\int_{\omega }\left\vert u\left( x,T\right) \right\vert ^{2}dx\right)
^{1-\mu }\left( \int_{\Omega }\left\vert u\left( x,0\right) \right\vert
^{2}dx\right) ^{\mu }\text{ .}
\end{equation*}%
Moreover, when $A$ is time-independent, then $c=Ce^{\frac{C}{T}}$ where $C$
and $\mu $ only depend on $\left( A,\Omega ,\omega ,n\right) $.
\end{thm}

\bigskip

\bigskip

Clearly, Proposition \ref{proposition2.1} is a direct application of Theorem %
\ref{theorem4.1}. The proof of Theorem \ref{theorem4.1} uses covering
argument and propagation of interpolation inequalities along a chain of
balls (also called propagation of smallness): First we extend $A\left( \cdot
,T\right) $ to a $C^{2}$ function on $\mathbb{R}^{n}$ denoted $A_{T}$. Next,
for each $x_{0}\in \mathbb{R}^{n}$ there are a neighborhood of $x_{0}$ and a
function $x\mapsto d\left( x,x_{0}\right) $ on which the following four
properties hold:

\begin{enumerate}
\item $\frac{1}{C}\left\vert x-x_{0}\right\vert \leq d\left( x,x_{0}\right)
\leq C\left\vert x-x_{0}\right\vert $ for some $C\geq1$ depending on $\left(
x_{0},A_{T}\right) $ ;

\item $x\mapsto d^{2}\left( x,x_{0}\right) $ is $C^{2}$ ;

\item $A_{T}\left( x\right) \nabla d\left( x,x_{0}\right) \cdot\nabla
d\left( x,x_{0}\right) =1$ ;

\item $\frac{1}{2}A_{T}\left( x\right) \nabla ^{2}d^{2}\left( x,x_{0}\right)
=I_{n}+O\left( d\left( x,x_{0}\right) \right) $ .
\end{enumerate}

Here $\nabla ^{2}$ denotes the Hessian matrix and $d\left( x,x_{0}\right) $
is the geodesic distance connecting $x$ to $x_{0}$. The proof of the above
properties for $d\left( x,x_{0}\right) $ is a consequence of Gauss's lemma
for $C^{2}$ metrics (see \cite[page 7]{IM}).

\bigskip

Now we are able to define the ball of center $x_{0}$ and radius $R$ as $%
B_{R}=\left\{ x;d\left( x,x_{0}\right) <R\right\} $. We will choose $%
x_{0}\in \Omega $ in order that one of the two following assumptions hold:
(i) $\overline{B_{r}}\subset \Omega $ for any $r\in \left( 0,R\right] $;
(ii) $B_{r}\cap \partial \Omega \neq \emptyset $ and $A\nabla d^{2}\cdot \nu
\geq 0$\textit{\ }on\textit{\ }$\partial \Omega \cap B_{R}$ for any $r\in %
\left[ R_{0},R\right] $ where $R_{0}>0$. Here $\nu $ is the unit outward
normal vector to $\partial \Omega \cap B_{R}$.

\bigskip

The case $\left( i\right) $ deals with the propagation in the interior
domain by a chain of balls strictly included in $\Omega $. The analysis near
the boundary $\partial \Omega $ requires the assumptions of $\left(
ii\right) $.

However when $\Omega \subset \mathbb{R}^{n}$ is a convex domain or a
star-shaped domain with respect to $x_{0}\in \Omega $, we only need to
propagate the estimate in the interior domain.

If further $A=I_{n}$, then $d\left( x,x_{0}\right) =\left\vert
x-x_{0}\right\vert $ and it is well defined for any $x\in \Omega $. From 
\cite{PW} such observation at one point in time implies the observability
for the heat equation which from \cite{AEWZ} is equivalent to the sum of
eigenfunctions estimate of Lebeau-Robbiano type. Eventually a careful
evaluation of the constants gives the following estimates (whose proof is
omitted).

\bigskip

\begin{thm}
\label{theorem4.2} Suppose that $\Omega \subset \mathbb{R}^{n}$ is a convex
domain or a star-shaped domain with respect to $x_{0}\in \Omega $ such that $%
\left\{ x;\left\vert x-x_{0}\right\vert <r\right\} \Subset \Omega $ for some 
$r>0$. Then for any $u_{0}\in L^{2}\left( \Omega \right) $,\thinspace $T>0$, 
$\left( a_{i}\right) _{i\geq 1}\in \mathbb{R}$, $\mu \geq 1$, $\varepsilon
\in \left( 0,1\right) $, one has%
\begin{equation*}
\left\Vert e^{T\Delta }u_{0}\right\Vert _{L^{2}\left( \Omega \right) }\leq 
\frac{1}{r^{n}}\frac{1}{r^{\varepsilon \left( n-2\right) }}e^{\frac{C}{T}%
\frac{1}{r^{6\varepsilon }}}\int_{0}^{T}\left\Vert e^{t\Delta
}u_{0}\right\Vert _{L^{2}\left( \left\vert x-x_{0}\right\vert <r\right) }dt%
\text{ }
\end{equation*}%
and%
\begin{equation*}
\sum\limits_{\mu _{i}\leq \mu }\left\vert a_{i}\right\vert ^{2}\leq \frac{1}{%
r^{2n\left( 1+\varepsilon \right) }}e^{C\frac{1}{r^{2\varepsilon }}\sqrt{\mu 
}}{\int\nolimits_{\left\vert x-x_{0}\right\vert <r}}\left\vert
\sum\limits_{\mu _{i}\leq \mu }a_{i}e_{i}\left( x\right) \right\vert ^{2}dx
\end{equation*}%
where $C>0$ is a constant only depending on $\left( \varepsilon ,n,%
\normalfont{\text{max}}\left\{ \left\vert x-x_{0}\right\vert ;x\in \overline{%
\Omega }\right\} \right) $. Here $\left( \mu _{i},e_{i}\right) $ denotes the
eigenbasis of the Laplace operator with Dirichlet boundary condition.
\end{thm}

\bigskip

\bigskip

In the next subsection, we state some preliminary lemmas and corollaries. In
subsection 4.2, we prove Theorem \ref{theorem4.1}. Subsection 4.3 is devoted
to the proof of the preliminary results.

\bigskip

\subsection{Preliminary results}

\bigskip

In this subsection we present some lemmas and corollaries which will be used
for the proof of Theorem \ref{theorem4.1}.

\bigskip

The following lemma allows to solve differential inequalities and makes
appear the H\"{o}lder type of inequality in Theorem \ref{theorem4.1}.

\bigskip

\begin{lem}
\label{lemma4.3} Let $T>0$, $\lambda >0$\ and $F_{1},F_{2}\in C^{0}\left( %
\left[ 0,T\right] \right) $. Consider two positive functions $y,N\in
C^{1}\left( \left[ 0,T\right] \right) $\ such that 
\begin{equation*}
\left\{ 
\begin{array}{ll}
\left\vert \displaystyle\frac{1}{2}y^{\prime }\left( t\right) +N\left(
t\right) y\left( t\right) \right\vert \leq \left( \displaystyle\frac{C_{0}}{%
T-t+\lambda }+C_{1}\right) y\left( t\right) +F_{1}\left( t\right) y\left(
t\right) &  \\ 
N^{\prime }\left( t\right) \leq \left( \displaystyle\frac{1+C_{0}}{%
T-t+\lambda }+C_{1}\right) N\left( t\right) +F_{2}\left( t\right) \text{ } & 
\end{array}%
\right.
\end{equation*}%
where $C_{0},C_{1}\geq 0$. Then for any $0\leq t_{1}<t_{2}<t_{3}\leq T$, one
has 
\begin{equation*}
y\left( t_{2}\right) ^{1+M}\leq y\left( t_{3}\right) y\left( t_{1}\right)
^{M}e^{4D}\left( \frac{T-t_{1}+\lambda }{T-t_{3}+\lambda }\right)
^{2C_{0}\left( 1+M\right) }
\end{equation*}%
where%
\begin{equation*}
M=\frac{\displaystyle\int_{t_{2}}^{t_{3}}\frac{e^{tC_{1}}}{\left(
T-t+\lambda \right) ^{1+C_{0}}}dt}{\displaystyle\int_{t_{1}}^{t_{2}}\frac{%
e^{tC_{1}}}{\left( T-t+\lambda \right) ^{1+C_{0}}}dt}\text{ and }D=M\left(
t_{2}-t_{1}\right) \left( C_{1}+\underset{\left[ t_{1},t_{3}\right] }{%
\normalfont{\text{sup}}}\left\vert F_{1}\right\vert +\displaystyle%
\int_{t_{1}}^{t_{3}}\left\vert F_{2}\right\vert dt\right) \text{ .}
\end{equation*}
\end{lem}

\bigskip

\begin{cor}
\label{corollary4.4} Under the assumptions of Lemma 1, for any $\lambda >0$\
and $\ell >1$\ such that $\ell \lambda <T/4$, one has 
\begin{equation*}
y\left( T-\ell \lambda \right) ^{1+M_{\ell }}\leq y\left( T\right) y\left(
T-2\ell \lambda \right) ^{M_{\ell }}e^{D_{\ell }}\left( 2\ell +1\right)
^{2C_{0}\left( 1+M_{\ell }\right) }\text{ }
\end{equation*}%
where $D_{\ell }=TM_{\ell }\left( C_{1}+\underset{\left[ t_{1},t_{3}\right] }%
{\normalfont{\text{sup}}}\left\vert F_{1}\right\vert +\displaystyle%
\int_{t_{1}}^{t_{3}}\left\vert F_{2}\right\vert dt\right) $, $M_{\ell }\leq
e^{C_{1}T}\frac{\left( \ell +1\right) ^{C_{0}}}{1-\left( \frac{2}{3}\right)
^{C_{0}}}$\ if $C_{0}>0$\ and $M_{\ell }\leq e^{C_{1}T}\frac{%
\normalfont{\text{ln}}\left( \ell +1\right) }{\normalfont{\text{ln}}2}$\ if $%
C_{0}=0$.
\end{cor}

\bigskip

Proof .- Apply Lemma \ref{lemma4.3} with $t_{3}=T$, $t_{2}=T-\ell \lambda $, 
$t_{1}=T-2\ell \lambda $, with $\ell \lambda <T/4$. Here when $C_{0}>0$ 
\begin{equation*}
M_{\ell }=\frac{\displaystyle\int_{T-\ell \lambda }^{T}\frac{e^{tC_{1}}}{%
\left( T-t+\lambda \right) ^{1+C_{0}}}dt}{\displaystyle\int_{T-2\ell \lambda
}^{T-\ell \lambda }\frac{e^{tC_{1}}}{\left( T-t+\lambda \right) ^{1+C_{0}}}dt%
}\leq e^{2\ell \lambda C_{1}}\frac{\left( \ell +1\right) ^{C_{0}}-1}{%
1-\left( \frac{\ell +1}{2\ell +1}\right) ^{C_{0}}}\leq e^{C_{1}T}\frac{%
\left( \ell +1\right) ^{C_{0}}}{1-\left( \frac{2}{3}\right) ^{C_{0}}}\text{
for }\ell >1\text{ .}
\end{equation*}%
And when $C_{0}=0$ 
\begin{equation*}
M_{\ell }=\frac{\displaystyle\int_{T-\ell \lambda }^{T}\frac{e^{tC_{1}}}{%
\left( T-t+\lambda \right) }dt}{\displaystyle\int_{T-2\ell \lambda }^{T-\ell
\lambda }\frac{e^{tC_{1}}}{\left( T-t+\lambda \right) }dt}\leq e^{2\ell
\lambda C_{1}}\frac{\text{ln}\left( \ell +1\right) }{\text{ln}\left( \frac{%
2\ell +1}{\ell +1}\right) }\leq e^{C_{1}T}\frac{\text{ln}\left( \ell
+1\right) }{\text{ln}2}\text{ for }\ell >1\text{ .}
\end{equation*}

\bigskip

The following lemma establishes the differential inequalities associated to
parabolic equations in any open set $\vartheta \subset \mathbb{R}^{n}$:

\bigskip

\begin{lem}
\label{lemma4.5} For any $\xi \in C^{2}\left( \overline{\Omega }\times \left[
0,T\right] \right) $, $z\in H^{1}\left( 0,T;H_{0}^{1}\left( \vartheta
\right) \right) $, one has%
\begin{equation*}
\begin{array}{ll}
& \quad \displaystyle\frac{1}{2}\frac{d}{dt}\int_{\vartheta }\left\vert
z\right\vert ^{2}e^{\xi }dx+\displaystyle\int_{\vartheta }A\nabla z\cdot
\nabla ze^{\xi }dx \\ 
& =\displaystyle\frac{1}{2}\int_{\vartheta }\left\vert z\right\vert
^{2}\left( \partial _{t}\xi +\nabla \cdot \left( A\nabla \xi \right)
+A\nabla \xi \cdot \nabla \xi \right) e^{\xi }dx+\displaystyle%
\int_{\vartheta }z\left( \partial _{t}z-\nabla \cdot \left( A\nabla z\right)
\right) e^{\xi }dx%
\end{array}%
\end{equation*}%
and for some $C$\ only depending on $\left( A,\partial _{x}A,\partial
_{t}A\right) $%
\begin{equation*}
\begin{array}{ll}
\displaystyle\frac{d}{dt}\frac{\displaystyle\int_{\vartheta }A\nabla z\cdot
\nabla ze^{\xi }dx}{\displaystyle\int_{\vartheta }\left\vert z\right\vert
^{2}e^{\xi }dx} & \leq \frac{-\displaystyle2\int_{\vartheta }A\nabla ^{2}\xi
A\nabla z\cdot \nabla ze^{\xi }dx}{\displaystyle\int_{\vartheta }\left\vert
z\right\vert ^{2}e^{\xi }dx}+\frac{\displaystyle\int_{\partial \vartheta
}\left( A\nabla z\cdot \nabla z\right) \left( A\nabla \xi \cdot \nu \right)
e^{\xi }dx}{\displaystyle\int_{\vartheta }\left\vert z\right\vert ^{2}e^{\xi
}dx} \\ 
& \quad +\frac{\displaystyle\int_{\vartheta }\left\vert \partial
_{t}z-\nabla \cdot \left( A\nabla z\right) \right\vert ^{2}e^{\xi }dx}{%
\displaystyle\int_{\vartheta }\left\vert z\right\vert ^{2}e^{\xi }dx}+C\frac{%
\displaystyle\int_{\vartheta }\left( 1+\left\vert \nabla \xi \right\vert
\right) \left\vert \nabla z\right\vert ^{2}e^{\xi }dx}{\displaystyle%
\int_{\vartheta }\left\vert z\right\vert ^{2}e^{\xi }dx} \\ 
& \quad +\frac{\displaystyle\int_{\vartheta }A\nabla z\cdot \nabla z\left(
\partial _{t}\xi +\nabla \cdot \left( A\nabla \xi \right) +A\nabla \xi \cdot
\nabla \xi \right) e^{\xi }dx}{\displaystyle\int_{\vartheta }\left\vert
z\right\vert ^{2}e^{\xi }dx} \\ 
& \quad -\frac{\displaystyle\int_{\vartheta }A\nabla z\cdot \nabla ze^{\xi
}dx}{\displaystyle\int_{\vartheta }\left\vert z\right\vert ^{2}e^{\xi }dx}%
\times \frac{\displaystyle\int_{\vartheta }\left\vert z\right\vert
^{2}\left( \partial _{t}\xi +\nabla \cdot \left( A\nabla \xi \right)
+A\nabla \xi \cdot \nabla \xi \right) e^{\xi }dx}{\displaystyle%
\int_{\vartheta }\left\vert z\right\vert ^{2}e^{\xi }dx}\text{ .}%
\end{array}%
\end{equation*}
\end{lem}

\bigskip

\begin{cor}
\label{corollary4.6} Let $R>0$\ be sufficiently small and $z\in H^{1}\left(
0,T;H_{0}^{1}\left( \Omega \cap B_{R}\right) \right) $\ with $B_{R}=\left\{
x;d\left( x,x_{0}\right) <R\right\} $. Introduce for $t\in \left( 0,T\right] 
$, $\mathcal{P}z=\partial _{t}z-\nabla \cdot \left( A\nabla z\right) $,%
\begin{equation*}
G_{\lambda }\left( x,t\right) =\frac{1}{\left( T-t+\lambda \right) ^{n/2}}%
e^{-\frac{d^{2}\left( x,x_{0}\right) }{4\left( T-t+\lambda \right) }}\quad
\forall x\in B_{R}\text{ ,}
\end{equation*}%
and 
\begin{equation*}
N_{\lambda }\left( t\right) =\frac{\displaystyle\int_{\Omega \cap
B_{R}}A\left( x,t\right) \nabla z\left( x,t\right) \cdot \nabla z\left(
x,t\right) G_{\lambda }\left( x,t\right) dx}{\displaystyle\int_{\Omega \cap
B_{R}}\left\vert z\left( x,t\right) \right\vert ^{2}G_{\lambda }\left(
x,t\right) dx}
\end{equation*}%
whenever $\displaystyle\int_{\Omega \cap B_{R}}\left\vert z\left( x,t\right)
\right\vert ^{2}dx\neq 0$. Then, the following two properties hold:

\begin{description}
\item[$i)$] For some $C_{0}\geq 0$,%
\begin{equation*}
\begin{array}{ll}
& \quad \left\vert \displaystyle\frac{1}{2}\frac{d}{dt}\int_{\Omega \cap
B_{R}}\left\vert z\left( x,t\right) \right\vert ^{2}G_{\lambda }\left(
x,t\right) dx+N_{\lambda }\left( t\right) \displaystyle\int_{\Omega \cap
B_{R}}\left\vert z\left( x,t\right) \right\vert ^{2}G_{\lambda }\left(
x,t\right) dx\right\vert \\ 
& \leq \left( \displaystyle\frac{C_{0}}{T-t+\lambda }+C_{1}\right)
\int_{\Omega \cap B_{R}}\left\vert z\left( x,t\right) \right\vert
^{2}G_{\lambda }\left( x,t\right) dx \\ 
& \quad +\displaystyle\int_{\Omega \cap B_{R}}\left\vert z\left( x,t\right) 
\mathcal{P}z\left( x,t\right) \right\vert G_{\lambda }\left( x,t\right) dx%
\text{ .}%
\end{array}%
\end{equation*}

\item[$ii)$] There are $R>0$, $0\leq C_{0}<1$, $C_{1}\geq 0$\ such that when 
$A\nabla d^{2}\cdot \nu \geq 0$\ on $\partial \Omega \cap B_{R}$,\ 
\begin{equation*}
\frac{d}{dt}N_{\lambda }\left( t\right) \leq \left( \frac{1+C_{0}}{%
T-t+\lambda }+C_{1}\right) N_{\lambda }\left( t\right) +\frac{\displaystyle%
\int_{\Omega \cap B_{R}}\left\vert \mathcal{P}z\left( x,t\right) \right\vert
^{2}G_{\lambda }\left( x,t\right) dx}{\displaystyle\int_{\Omega \cap
B_{R}}\left\vert z\left( x,t\right) \right\vert ^{2}G_{\lambda }\left(
x,t\right) dx}\text{ .}
\end{equation*}
\end{description}
\end{cor}

\bigskip

Proof .- Apply Lemma \ref{lemma4.5} with $\vartheta =\Omega \cap B_{R}$ and 
\begin{equation*}
\xi \left( x,t\right) =-\frac{d^{2}\left( x,x_{0}\right) }{4\left(
T-t+\lambda \right) }-\frac{n}{2}\text{ln}\left( T-t+\lambda \right) \text{ .%
}
\end{equation*}%
It remains to bound $\partial _{t}\xi +\nabla \cdot \left( A\nabla \xi
\right) +A\nabla \xi \cdot \nabla \xi $, $-2A\nabla ^{2}\xi A\nabla z\cdot
\nabla z$ and $\left\vert \nabla \xi \right\vert $. First, one get 
\begin{equation*}
C_{A_{T}}\left\vert \nabla \xi \right\vert ^{2}\leq A_{T}\nabla \xi \cdot
\nabla \xi =\frac{d^{2}\left( x,x_{0}\right) }{4\left( T-t+\lambda \right)
^{2}}\text{ .}
\end{equation*}%
Next, $\nabla \cdot \left( A_{T}\nabla \xi \right) =\frac{-n}{\left(
T-t+\lambda \right) }+\frac{O\left( d\left( x,x_{0}\right) \right) }{%
T-t+\lambda }$ and 
\begin{equation*}
\partial _{t}\xi +\nabla \cdot \left( A_{T}\nabla \xi \right) +A_{T}\nabla
\xi \cdot \nabla \xi =\frac{O\left( d\left( x,x_{0}\right) \right) }{%
T-t+\lambda }\text{ }
\end{equation*}%
imply 
\begin{equation*}
\begin{array}{ll}
& \quad \partial _{t}\xi +\nabla \cdot \left( A\nabla \xi \right) +A\nabla
\xi \cdot \nabla \xi \\ 
& =\frac{O\left( d\left( x,x_{0}\right) \right) }{T-t+\lambda }+\nabla \cdot
\left( \left( A-A_{T}\right) \left( A_{T}\right) ^{-1}A_{T}\nabla \xi
\right) +\left( A-A_{T}\right) \nabla \xi \cdot \nabla \xi \\ 
& =\frac{O\left( d\left( x,x_{0}\right) \right) }{T-t+\lambda }+O\left(
1\right) \text{ ,}%
\end{array}%
\end{equation*}%
where in the last equality we used $\left\Vert A\left( \cdot ,t\right)
-A_{T}\right\Vert \leq \left\Vert \partial _{t}A\right\Vert \left(
T-t+\lambda \right) $. Finally, we have%
\begin{equation*}
\begin{array}{ll}
-2A\nabla ^{2}\xi A\nabla z\cdot \nabla z & =\frac{1}{2\left( T-t+\lambda
\right) }A_{T}\nabla ^{2}d^{2}A\nabla z\cdot \nabla z+\left( A-A_{T}\right)
\nabla ^{2}d^{2}A\nabla z\cdot \nabla z \\ 
& =A\nabla z\cdot \nabla z\left( \frac{1+O\left( d\left( x,x_{0}\right)
\right) }{T-t+\lambda }+O\left( 1\right) \right) \text{ .}%
\end{array}%
\end{equation*}%
One conclude by choosing $R>0$ sufficiently small in order the constant $%
C_{0}$ in Corollary \ref{corollary4.6} satisfies $0<C_{0}<1$.

\bigskip

Remark .- When $A$ is time-independent, then $C_{1}=0$ in Corollary \ref%
{corollary4.6}.

\bigskip

The following lemma will be used to deal with the delocalized terms.

\bigskip

\begin{lem}
\label{lemma4.7} Let $\rho \in \left( 0,R\right) $\ and $0<\varepsilon <\rho
/2$. There are constants $c_{1},c_{2},c_{3}>0$ only depending on $\left(
\rho ,\varepsilon ,A\right) $\ such that for any $T-\theta \leq t\leq T$,
one has 
\begin{equation*}
\frac{\displaystyle\int_{\Omega }\left\vert u\left( x,0\right) \right\vert
^{2}dx}{\displaystyle\int_{\Omega \cap B_{\rho }}\left\vert u\left(
x,t\right) \right\vert ^{2}dx}\leq e^{\frac{c_{1}}{\theta }}\text{ }
\end{equation*}%
where 
\begin{equation*}
\frac{1}{\theta }=c_{2}\normalfont{\text{ln}}\left( e^{c_{3}\left( 1+\frac{1%
}{T}\right) }\frac{\displaystyle\int_{\Omega }\left\vert u\left( x,0\right)
\right\vert ^{2}dx}{\displaystyle\int_{\Omega \cap B_{\rho -2\varepsilon
}}\left\vert u\left( x,T\right) \right\vert ^{2}dx}\right) \text{\textit{\
with }}0<\theta \leq \normalfont{\text{min}}\left( 1,T/2\right) \text{ .}
\end{equation*}
\end{lem}

\bigskip

The interested reader may wish here to compare this lemma with \cite[Lemma 5]%
{EFV}.

\bigskip

\bigskip

\subsection{Proof of Theorem 4.1}

\bigskip

Let $\lambda >0$ and $\ell >1$ be such that $\ell \lambda <T/4$. By
Corollary \ref{corollary4.4} with $y\left( t\right) =\displaystyle%
\int_{\Omega \cap B_{R}}\left\vert z\left( x,t\right) \right\vert
^{2}G_{\lambda }\left( x,t\right) dx$, $N\left( t\right) =N_{\lambda }\left(
t\right) $ given in Corollary \ref{corollary4.6}, 
\begin{equation*}
F_{1}\left( t\right) =\frac{\displaystyle\int_{\Omega \cap B_{R}}\left\vert
z\left( x,t\right) \left( \partial _{t}z\left( x,t\right) -\nabla \cdot
\left( A\left( x,t\right) \nabla z\left( x,t\right) \right) \right)
\right\vert G_{\lambda }\left( x,t\right) dx}{\displaystyle\int_{\Omega \cap
B_{R}}\left\vert z\left( x,t\right) \right\vert ^{2}G_{\lambda }\left(
x,t\right) dx}
\end{equation*}%
and 
\begin{equation*}
F_{2}\left( t\right) =\frac{\displaystyle\int_{\Omega \cap B_{R}}\left\vert
\partial _{t}z\left( x,t\right) -\nabla \cdot \left( A\left( x,t\right)
\nabla z\left( x,t\right) \right) \right\vert ^{2}G_{\lambda }\left(
x,t\right) dx}{\displaystyle\int_{\Omega \cap B_{R}}\left\vert z\left(
x,t\right) \right\vert ^{2}G_{\lambda }\left( x,t\right) dx}\text{ }
\end{equation*}%
knowing that $N^{\prime }\left( t\right) \leq \left( \frac{1+C_{0}}{%
T-t+\lambda }+C_{1}\right) N\left( t\right) +F_{2}\left( t\right) $ from
Corollary \ref{corollary4.6}, one can deduce the following interpolation
inequality with $M_{\ell }\leq e^{C_{1}T}\frac{\left( \ell +1\right) ^{C_{0}}%
}{1-\left( \frac{2}{3}\right) ^{C_{0}}}$ and $0<C_{0}<1$,%
\begin{equation*}
y\left( T-\ell \lambda \right) ^{1+M_{\ell }}\leq y\left( T\right) y\left(
T-2\ell \lambda \right) ^{M_{\ell }}\left( 2\ell +1\right) ^{2C_{0}\left(
1+M_{\ell }\right) }e^{D_{\ell }}\text{ }
\end{equation*}%
that is 
\begin{equation*}
\begin{array}{ll}
& \quad \left( \displaystyle\int_{\Omega \cap B_{R}}\left\vert z\left(
x,T-\ell \lambda \right) \right\vert ^{2}e^{\frac{-d^{2}\left(
x,x_{0}\right) }{4\left( \ell +1\right) \lambda }}dx\right) ^{1+M_{\ell }}
\\ 
& \leq \left( \ell +1\right) ^{n/2}\left( 2\ell +1\right) ^{2C_{0}\left(
1+M_{\ell }\right) }\displaystyle\int_{\Omega \cap B_{R}}\left\vert z\left(
x,T\right) \right\vert ^{2}e^{\frac{-d^{2}\left( x,x_{0}\right) }{4\lambda }%
}dx\left( \displaystyle\int_{\Omega }\left\vert u\left( x,0\right)
\right\vert ^{2}dx\right) ^{M_{\ell }} \\ 
& \quad \times e^{TM_{\ell }\left( \displaystyle\int_{T-2\ell \lambda
}^{T}\left\vert F_{2}\right\vert dt+\underset{\left[ T-2\ell \lambda ,T%
\right] }{\text{sup}}\left\vert F_{1}\right\vert +C_{1}\right) }\text{ .}%
\end{array}%
\end{equation*}%
From the definition of $F_{1}$, 
\begin{equation*}
\left\vert F_{1}\left( t\right) \right\vert \leq e^{\frac{\left(
R-2\varepsilon \right) ^{2}}{4\left( T-t+\lambda \right) }}e^{-\frac{\left(
R-\varepsilon \right) ^{2}}{4\left( T-t+\lambda \right) }}\frac{\displaystyle%
\int_{\Omega \cap \left\{ R-\varepsilon \leq d\left( x,x_{0}\right) \right\}
}\left\vert \chi u\right\vert \left\vert -2A\nabla \chi \cdot \nabla
u-\nabla \cdot \left( A\nabla \chi \right) u\right\vert dx}{\displaystyle%
\int_{\Omega \cap B_{R-2\varepsilon }}\left\vert u\right\vert ^{2}dx}\text{ .%
}
\end{equation*}%
Since $e^{\frac{\left( R-2\varepsilon \right) ^{2}}{4\left( T-t+\lambda
\right) }}e^{-\frac{\left( R-\varepsilon \right) ^{2}}{4\left( T-t+\lambda
\right) }}=e^{-\frac{\varepsilon \left( 2R-3\varepsilon \right) }{4\left(
T-t+\lambda \right) }}\leq e^{-\frac{\varepsilon \left( 2R-3\varepsilon
\right) }{12\ell \lambda }}$ for $t\in \left[ T-2\ell \lambda ,T\right] $
with $\ell >1$, one has when $t\in \left[ T-2\ell \lambda ,T\right] $ 
\begin{equation*}
\left\vert F_{1}\left( t\right) \right\vert \leq e^{-\frac{\varepsilon
\left( 2R-3\varepsilon \right) }{12\ell \lambda }}\frac{\displaystyle%
\int_{\Omega \cap \left\{ R-\varepsilon \leq d\left( x,x_{0}\right) \right\}
}\left\vert \chi u\right\vert \left\vert -2A\nabla \chi \cdot \nabla
u-\nabla \cdot \left( A\nabla \chi \right) u\right\vert dx}{\displaystyle%
\int_{\Omega \cap B_{R-2\varepsilon }}\left\vert u\right\vert ^{2}dx}\text{ .%
}
\end{equation*}%
By Lemma \ref{lemma4.7} with $\rho =R-2\varepsilon $,%
\begin{equation*}
\underset{t\in \left[ T-2\ell \lambda ,T\right] }{\text{sup}}\left\vert
F_{1}\left( t\right) \right\vert \leq e^{-\frac{\varepsilon \left(
2R-3\varepsilon \right) }{12\ell \lambda }}ce^{\frac{c_{1}}{\theta }}\text{
if }2\ell \lambda \leq \theta \text{\ .}
\end{equation*}%
Similarly, from the definition of $F_{2}$, 
\begin{equation*}
\left\vert F_{2}\left( t\right) \right\vert \leq e^{\frac{\left(
R-2\varepsilon \right) ^{2}}{4\left( T-t+\lambda \right) }}e^{-\frac{\left(
R-\varepsilon \right) ^{2}}{4\left( T-t+\lambda \right) }}\frac{\displaystyle%
\int_{\Omega \cap \left\{ R-\varepsilon \leq d\left( x,x_{0}\right) \right\}
}\left\vert -2A\nabla \chi \cdot \nabla u-\nabla \cdot \left( A\nabla \chi
\right) u\right\vert ^{2}dx}{\displaystyle\int_{\Omega \cap
B_{R-2\varepsilon }}\left\vert u\right\vert ^{2}dx}
\end{equation*}%
and then, when $t\in \left[ T-2\ell \lambda ,T\right] $ 
\begin{equation*}
\left\vert F_{2}\left( t\right) \right\vert \leq e^{-\frac{\varepsilon
\left( 2R-3\varepsilon \right) }{12\ell \lambda }}\frac{\displaystyle%
\int_{\Omega \cap \left\{ R-\varepsilon \leq d\left( x,x_{0}\right) \right\}
}\left\vert -2A\nabla \chi \cdot \nabla u-\nabla \cdot \left( A\nabla \chi
\right) u\right\vert ^{2}dx}{\displaystyle\int_{\Omega \cap
B_{R-2\varepsilon }}\left\vert u\right\vert ^{2}dx}\text{ .}
\end{equation*}%
By Lemma \ref{lemma4.7} with $\rho =R-2\varepsilon $,%
\begin{equation*}
\int_{T-2\ell \lambda }^{T}\left\vert F_{2}\left( t\right) \right\vert
dt\leq e^{-\frac{\varepsilon \left( 2R-3\varepsilon \right) }{12\ell \lambda 
}}ce^{\frac{c_{1}}{\theta }}\text{ if }2\ell \lambda \leq \theta \text{ }
\end{equation*}%
where $c>1$ is a constant only dependent on $\left( A,R,\varepsilon \right) $%
. We conclude that for any $2\ell \lambda \leq \theta \frac{\varepsilon
\left( 2R-3\varepsilon \right) }{6c_{1}}$%
\begin{equation*}
\underset{t\in \left[ T-\theta ,T\right] }{\text{sup}}\left\vert F_{1}\left(
t\right) \right\vert +\int_{T-\theta }^{T}\left\vert F_{2}\left( t\right)
\right\vert dt\leq 2c\text{ .}
\end{equation*}%
Therefore there is $c_{4}:=\frac{\varepsilon \left( 2R-3\varepsilon \right) 
}{6c_{1}}\in \left( 0,1\right) $ such that for any $2\ell \lambda \leq
c_{4}\theta $ 
\begin{equation*}
\begin{array}{ll}
& \quad \left( \displaystyle\int_{\Omega \cap B_{R}}\left\vert z\left(
x,T-\ell \lambda \right) \right\vert ^{2}e^{\frac{-d^{2}\left(
x,x_{0}\right) }{4\left( \ell +1\right) \lambda }}dx\right) ^{1+M_{\ell }}
\\ 
& \leq e^{\left( 2c+C_{1}\right) TM_{\ell }}\left( 2\ell +1\right)
^{2C_{0}\left( 1+M_{\ell }\right) +n/2}\displaystyle\int_{\Omega \cap
B_{R}}\left\vert z\left( x,T\right) \right\vert ^{2}e^{\frac{-d^{2}\left(
x,x_{0}\right) }{4\lambda }}dx\left( \displaystyle\int_{\Omega }\left\vert
u\left( x,0\right) \right\vert ^{2}dx\right) ^{M_{\ell }}%
\end{array}%
\end{equation*}%
which implies%
\begin{equation*}
\begin{array}{ll}
\left( \displaystyle\int_{\Omega \cap B_{R}}\left\vert z\left( x,T-\ell
\lambda \right) \right\vert ^{2}dx\right) ^{1+M_{\ell }} & \leq e^{\left(
2c+C_{1}\right) TM_{\ell }}\left( 2\ell +1\right) ^{2C_{0}\left( 1+M_{\ell
}\right) +n/2}e^{\frac{R^{2}}{4\left( \ell +1\right) \lambda }\left(
1+M_{\ell }\right) } \\ 
& \quad \times \displaystyle\int_{\Omega \cap B_{R}}\left\vert z\left(
x,T\right) \right\vert ^{2}e^{\frac{-d^{2}\left( x,x_{0}\right) }{4\lambda }%
}dx\left( \displaystyle\int_{\Omega }\left\vert u\left( x,0\right)
\right\vert ^{2}dx\right) ^{M_{\ell }}\text{ .}%
\end{array}%
\end{equation*}%
Now, we split $\displaystyle\int_{\Omega \cap B_{R}}\left\vert z\left(
x,T\right) \right\vert ^{2}e^{\frac{-d^{2}\left( x,x_{0}\right) }{4\lambda }%
}dx$ into two parts: For any $0<r<R/2$ such that $B_{r}\Subset \Omega $, 
\begin{equation*}
\int_{\Omega \cap B_{R}}\left\vert z\left( x,T\right) \right\vert ^{2}e^{%
\frac{-d^{2}\left( x,x_{0}\right) }{4\lambda }}dx\leq \int_{B_{r}}\left\vert
u\left( x,T\right) \right\vert ^{2}dx+e^{\frac{-r^{2}}{4\lambda }%
}\int_{\Omega }\left\vert u\left( x,0\right) \right\vert ^{2}dx\text{ .}
\end{equation*}%
Consequently, we have%
\begin{equation*}
\begin{array}{ll}
& \quad \left( \displaystyle\int_{\Omega \cap B_{R}}\left\vert z\left(
x,T-\ell \lambda \right) \right\vert ^{2}dx\right) ^{1+M_{\ell }} \\ 
& \leq e^{\left( 2c+C_{1}\right) TM_{\ell }}\left( 2\ell +1\right)
^{2C_{0}\left( 1+M_{\ell }\right) +n/2}e^{\frac{R^{2}}{4\left( \ell
+1\right) \lambda }\left( 1+M_{\ell }\right) }\displaystyle%
\int_{B_{r}}\left\vert u\left( x,T\right) \right\vert ^{2}dx\left( %
\displaystyle\int_{\Omega }\left\vert u\left( x,0\right) \right\vert
^{2}dx\right) ^{M_{\ell }} \\ 
& \quad +e^{\left( 2c+C_{1}\right) TM_{\ell }}\left( 2\ell +1\right)
^{2C_{0}\left( 1+M_{\ell }\right) +n/2}e^{\frac{R^{2}}{4\left( \ell
+1\right) \lambda }\left( 1+M_{\ell }\right) }e^{\frac{-r^{2}}{4\lambda }%
}\left( \displaystyle\int_{\Omega }\left\vert u\left( x,0\right) \right\vert
^{2}dx\right) ^{1+M_{\ell }}\text{ .}%
\end{array}%
\end{equation*}%
Now, choose $\ell >1$ in order that $\frac{R^{2}}{4\left( \ell +1\right) }%
\left( 1+M_{\ell }\right) \leq \frac{r^{2}}{8}$ (knowing that $M_{\ell }\leq
e^{C_{1}T}\frac{\left( \ell +1\right) ^{C_{0}}}{1-\left( \frac{2}{3}\right)
^{C_{0}}}$ for $\ell >1$ and $C_{0}<1$). Therefore there is $K>1$ such that
for any $\lambda \leq \frac{c_{4}}{2\ell }\theta $ 
\begin{equation*}
\begin{array}{ll}
\left( \displaystyle\int_{\Omega \cap B_{R-\varepsilon }}\left\vert u\left(
x,T-\ell \lambda \right) \right\vert ^{2}dx\right) ^{1+K} & \leq Ke^{\frac{%
r^{2}}{8\lambda }}\displaystyle\int_{B_{r}}\left\vert u\left( x,T\right)
\right\vert ^{2}dx\left( \displaystyle\int_{\Omega }\left\vert u\left(
x,0\right) \right\vert ^{2}dx\right) ^{K} \\ 
& \quad +Ke^{\frac{-r^{2}}{8\lambda }}\left( \displaystyle\int_{\Omega
}\left\vert u\left( x,0\right) \right\vert ^{2}dx\right) ^{1+K}\text{ .}%
\end{array}%
\end{equation*}%
But by Lemma \ref{lemma4.7} with $\rho =R-2\varepsilon $, since $\ell
\lambda \leq \theta $, 
\begin{equation*}
\int_{\Omega }\left\vert u\left( x,T\right) \right\vert ^{2}dx\leq
\int_{\Omega }\left\vert u\left( x,0\right) \right\vert ^{2}dx\leq e^{\frac{%
c_{1}}{\theta }}\int_{\Omega \cap B_{R-\varepsilon }}\left\vert u\left(
x,T-\ell \lambda \right) \right\vert ^{2}dx\text{ .}
\end{equation*}%
As a consequence, for any $\lambda \leq \frac{c_{4}}{2\ell }\theta $ one
obtain%
\begin{equation*}
\begin{array}{ll}
\left( \displaystyle\int_{\Omega }\left\vert u\left( x,T\right) \right\vert
^{2}dx\right) ^{1+K} & \leq e^{\frac{\left( 1+K\right) c_{1}}{\theta }}Ke^{%
\frac{r^{2}}{8\lambda }}\displaystyle\int_{B_{r}}\left\vert u\left(
x,T\right) \right\vert ^{2}dx\left( \displaystyle\int_{\Omega }\left\vert
u\left( x,0\right) \right\vert ^{2}dx\right) ^{K} \\ 
& \quad +e^{\frac{\left( 1+K\right) c_{1}}{\theta }}Ke^{\frac{-r^{2}}{%
8\lambda }}\left( \displaystyle\int_{\Omega }\left\vert u\left( x,0\right)
\right\vert ^{2}dx\right) ^{1+K}\text{ .}%
\end{array}%
\end{equation*}%
On the other hand, for any $\lambda \in \left( \frac{c_{4}}{2\ell }\theta ,%
\frac{T}{4\ell }\right) $, one has $1\leq e^{\frac{-r^{2}}{8\lambda }}e^{%
\frac{r^{2}\ell }{4c_{4}\theta }}$. And for any $\lambda \geq \frac{T}{4\ell 
}$, there holds $1\leq e^{\frac{-r^{2}}{8\lambda }}e^{\frac{r^{2}\ell }{4T}}$%
. Finally, there is $K>1$ such that for any $\lambda >0$,%
\begin{equation*}
\begin{array}{ll}
\left( \displaystyle\int_{\Omega }\left\vert u\left( x,T\right) \right\vert
^{2}dx\right) ^{1+K} & \leq e^{\frac{K}{\theta }}Ke^{\frac{r^{2}}{8\lambda }}%
\displaystyle\int_{B_{r}}\left\vert u\left( x,T\right) \right\vert
^{2}dx\left( \displaystyle\int_{\Omega }\left\vert u\left( x,0\right)
\right\vert ^{2}dx\right) ^{K} \\ 
& \quad +e^{\frac{K}{\theta }}Ke^{\frac{K}{T}}e^{\frac{-r^{2}}{8\lambda }%
}\left( \displaystyle\int_{\Omega }\left\vert u\left( x,0\right) \right\vert
^{2}dx\right) ^{1+K}\text{ .}%
\end{array}%
\end{equation*}%
Next, choose $\lambda >0$ such that $e^{\frac{r^{2}}{8\lambda }}:=2e^{\frac{K%
}{\theta }}Ke^{\frac{K}{T}}\left( \frac{\displaystyle\int_{\Omega
}\left\vert u\left( x,0\right) \right\vert ^{2}dx}{\displaystyle\int_{\Omega
}\left\vert u\left( x,T\right) \right\vert ^{2}dx}\right) ^{1+K}$ that is 
\begin{equation*}
e^{\frac{K}{\theta }}Ke^{\frac{K}{T}}e^{\frac{-r^{2}}{8\lambda }}\left( %
\displaystyle\int_{\Omega }\left\vert u\left( x,0\right) \right\vert
^{2}dx\right) ^{1+K}=\frac{1}{2}\left( \displaystyle\int_{\Omega }\left\vert
u\left( x,T\right) \right\vert ^{2}dx\right) ^{1+K}
\end{equation*}%
in order that%
\begin{equation*}
\displaystyle\int_{\Omega }\left\vert u\left( x,T\right) \right\vert
^{2}dx\leq 2Ke^{\frac{K}{\theta }}\left( e^{\frac{K}{T}}\displaystyle%
\int_{B_{r}}\left\vert u\left( x,T\right) \right\vert ^{2}dx\right) ^{\frac{1%
}{2+2K}}\left( \displaystyle\int_{\Omega }\left\vert u\left( x,0\right)
\right\vert ^{2}dx\right) ^{\frac{1+2K}{2+2K}}\text{ .}
\end{equation*}%
Recall that by Lemma \ref{lemma4.7} with $\rho =R-2\varepsilon $, 
\begin{equation*}
e^{\frac{K}{\theta }}=\left( e^{c_{3}\left( 1+\frac{1}{T}\right) }\frac{%
\displaystyle\int_{\Omega }\left\vert u\left( x,0\right) \right\vert ^{2}dx}{%
\displaystyle\int_{\Omega \cap B_{R-4\varepsilon }}\left\vert u\left(
x,T\right) \right\vert ^{2}dx}\right) ^{Kc_{2}}\text{ .}
\end{equation*}%
Finally, we obtain 
\begin{equation*}
\int_{\Omega \cap B_{R-4\varepsilon }}\left\vert u\left( x,T\right)
\right\vert ^{2}dx\leq K\left( \int_{\Omega }\left\vert u\left( x,0\right)
\right\vert ^{2}dx\right) ^{\frac{K}{1+K}}\left( e^{\frac{K}{T}%
}\int_{B_{r}}\left\vert u\left( x,T\right) \right\vert ^{2}dx\right) ^{\frac{%
1}{1+K}}
\end{equation*}%
for some positive constant $K$ only depending on $\left( A_{T},\varepsilon
,R,r,n\right) $. By an adequate covering of $\Omega $ by balls $%
B_{R-4\varepsilon }$ where $x_{0}$ and $R$ are chosen such that $A\nabla
d^{2}\cdot \nu \geq 0$\textit{\ }on\textit{\ }$\partial \Omega \cap B_{R}$
and by a propagation of smallness based on the previous estimate, we get the
desired observation inequality at one point in time for parabolic equations.

\bigskip

\bigskip

\subsection{Proof of Lemma \protect\ref{lemma4.3}}

\bigskip

We shall distinguish two cases: $t\in \left[ t_{1},t_{2}\right] $; $t\in %
\left[ t_{2},t_{3}\right] $. For $t_{1}\leq t\leq t_{2}$, we integrate $%
\left( \left( T-t+\lambda \right) ^{1+C_{0}}e^{-tC_{1}}N\left( t\right)
\right) ^{\prime }\leq \left( T-t+\lambda \right)
^{1+C_{0}}e^{-tC_{1}}F_{2}\left( t\right) $ over $\left( t,t_{2}\right) $ to
get 
\begin{equation*}
\left( \frac{T-t_{2}+\lambda }{T-t+\lambda }\right)
^{1+C_{0}}e^{-C_{1}\left( t_{2}-t\right) }N\left( t_{2}\right)
-\int_{t_{1}}^{t_{2}}\left\vert F_{2}\left( s\right) \right\vert ds\leq
N\left( t\right) \text{ .}
\end{equation*}%
Then we solve $y^{\prime }+2\alpha \left( t\right) y\leq 0$ with 
\begin{equation*}
\alpha \left( t\right) =\left( \frac{T-t_{2}+\lambda }{T-t+\lambda }\right)
^{1+C_{0}}e^{-C_{1}\left( t_{2}-t\right) }N\left( t_{2}\right) -\frac{C_{0}}{%
T-t+\lambda }-C_{1}-\int_{t_{1}}^{t_{2}}\left\vert F_{2}\right\vert ds-%
\underset{\left[ t_{1},t_{2}\right] }{\text{sup}}\left\vert F_{1}\right\vert
\end{equation*}%
and integrate it over $\left( t_{1},t_{2}\right) $ to obtain 
\begin{equation*}
\begin{array}{ll}
& \quad e^{\displaystyle2N\left( t_{2}\right) \int_{t_{1}}^{t_{2}}\left( 
\frac{T-t_{2}+\lambda }{T-t+\lambda }\right) ^{1+C_{0}}e^{-C_{1}\left(
t_{2}-t\right) }dt} \\ 
& \leq \frac{y\left( t_{1}\right) }{y\left( t_{2}\right) }\left( \frac{%
T-t_{1}+\lambda }{T-t_{2}+\lambda }\right) ^{2C_{0}}e^{2\left(
t_{2}-t_{1}\right) \left( \displaystyle C_{1}+\int_{t_{1}}^{t_{2}}\left\vert
F_{2}\right\vert ds+\underset{\left[ t_{1},t_{2}\right] }{\text{sup}}%
\left\vert F_{1}\right\vert \right) }\text{ .}%
\end{array}%
\end{equation*}%
For $t_{2}\leq t\leq t_{3}$, we integrate $\left( \left( T-t+\lambda \right)
^{1+C_{0}}e^{-tC_{1}}N\left( t\right) \right) ^{\prime }\leq \left(
T-t+\lambda \right) ^{1+C_{0}}F_{2}\left( t\right) $ over $\left(
t_{2},t\right) $ to get 
\begin{equation*}
N\left( t\right) \leq e^{C_{1}\left( t-t_{2}\right) }\left( \frac{%
T-t_{2}+\lambda }{T-t+\lambda }\right) ^{1+C_{0}}\left( N\left( t_{2}\right)
+\int_{t_{2}}^{t_{3}}\left\vert F_{2}\left( s\right) \right\vert ds\right) 
\text{ .}
\end{equation*}%
Then we solve $0\leq y^{\prime }+2\alpha \left( t\right) y$ with 
\begin{equation*}
\alpha \left( t\right) =\left( \frac{T-t_{2}+\lambda }{T-t+\lambda }\right)
^{1+C_{0}}e^{C_{1}\left( t-t_{2}\right) }(N\left( t_{2}\right)
+\int_{t_{2}}^{t_{3}}\left\vert F_{2}\right\vert ds+\underset{\left[
t_{2},t_{3}\right] }{\text{sup}}\left\vert F_{1}\right\vert +C_{1})+\frac{%
C_{0}}{T-t+\lambda }
\end{equation*}%
and integrate it over $\left( t_{2},t_{3}\right) $ to obtain 
\begin{equation*}
\begin{array}{ll}
y\left( t_{2}\right) & \leq e^{2\left( \displaystyle N\left( t_{2}\right)
+\int_{t_{2}}^{t_{3}}\left\vert F_{2}\right\vert ds+\underset{\left[
t_{2},t_{3}\right] }{\text{sup}}\left\vert F_{1}\right\vert +C_{1}\right) %
\displaystyle\int_{t_{2}}^{t_{3}}\left( \frac{T-t_{2}+\lambda }{T-t+\lambda }%
\right) ^{1+C_{0}}e^{C_{1}\left( t-t_{2}\right) }dt} \\ 
& \quad \times y\left( t_{3}\right) \left( \frac{T-t_{2}+\lambda }{%
T-t_{3}+\lambda }\right) ^{2C_{0}}\text{ .}%
\end{array}%
\end{equation*}%
Finally, combining the case $t_{1}\leq t\leq t_{2}$ and the case $t_{2}\leq
t\leq t_{3}$, we have 
\begin{equation*}
\begin{array}{ll}
y\left( t_{2}\right) & \leq y\left( t_{3}\right) \left( \frac{y\left(
t_{1}\right) }{y\left( t_{2}\right) }\right) ^{M}\left( \frac{%
T-t_{2}+\lambda }{T-t_{3}+\lambda }\right) ^{2C_{0}}\left( \frac{%
T-t_{1}+\lambda }{T-t_{2}+\lambda }\right) ^{2C_{0}M} \\ 
& \quad \times e^{\displaystyle2M\left( t_{2}-t_{1}\right)
\int_{t_{1}}^{t_{2}}\left\vert F_{2}\right\vert ds}e^{\displaystyle2M\left(
t_{2}-t_{1}\right) \int_{t_{2}}^{t_{3}}\left\vert F_{2}\right\vert ds} \\ 
& \quad \times e^{2M\left( t_{2}-t_{1}\right) \left( \underset{\left[
t_{1},t_{2}\right] }{\text{sup}}\left\vert F_{1}\right\vert +C_{1}\right)
}e^{2M\left( t_{2}-t_{1}\right) \left( \underset{\left[ t_{2},t_{3}\right] }{%
\text{sup}}\left\vert F_{1}\right\vert +C_{1}\right) }%
\end{array}%
\end{equation*}%
with%
\begin{equation*}
M=\frac{\displaystyle\int_{t_{2}}^{t_{3}}\frac{e^{tC_{1}}}{\left(
T-t+\lambda \right) ^{1+C_{0}}}dt}{\displaystyle\int_{t_{1}}^{t_{2}}\frac{%
e^{tC_{1}}}{\left( T-t+\lambda \right) ^{1+C_{0}}}dt}
\end{equation*}%
which is the desired estimate.

\bigskip

\bigskip

\subsection{Proof of Lemma \protect\ref{lemma4.5}}

\bigskip

The aim of this section is to prove the differential inequalities for
parabolic equations stated in Lemma \ref{lemma4.5}. For any $z\in
H^{1}\left( 0,T;H_{0}^{1}\left( \vartheta \right) \right) $, a weak solution
of $\partial _{t}z-\nabla \cdot \left( A\nabla z\right) =g$ with $g\in
L^{2}\left( \Omega \times \left( 0,T\right) \right) $, we apply the
following formula 
\begin{equation*}
\int_{\vartheta }{\partial }_{t}z\varphi dx+\int_{\vartheta }A\nabla z\cdot
\nabla \varphi dx=\int_{\partial \vartheta }A\nabla z\cdot \nu \varphi
dx+\int_{\vartheta }g\varphi dx\text{ }
\end{equation*}%
with different functions $\varphi $: $\varphi =ze^{\xi }$, $\varphi ={%
\partial }_{t}ze^{\xi }$ and $\varphi =A\nabla z\cdot \nabla \xi e^{\xi }$.
Here $\nu $ is the unit outward normal vector to $\partial \vartheta $ and $%
\xi =\xi \left( x,t\right) $ is a sufficiently smooth function which will be
chosen later. When $\varphi =ze^{\xi }$, we have 
\begin{equation*}
\int_{\vartheta }A\nabla z\cdot \nabla ze^{\xi }dx=-\int_{\vartheta }\left( {%
\partial }_{t}z+A\nabla z\cdot \nabla \xi -\frac{1}{2}g\right) ze^{\xi }dx+%
\frac{1}{2}\int_{\vartheta }gze^{\xi }dx
\end{equation*}%
and 
\begin{equation*}
\begin{array}{ll}
& \quad \displaystyle\frac{1}{2}\frac{d}{dt}\int_{\vartheta }\left\vert
z\right\vert ^{2}e^{\xi }dx=\displaystyle\int_{\vartheta }z{\partial }%
_{t}ze^{\xi }dx+\displaystyle\frac{1}{2}\int_{\vartheta }\left\vert
z\right\vert ^{2}\partial _{t}\xi e^{\xi }dx \\ 
& =-\displaystyle\int_{\vartheta }A\nabla z\cdot \nabla ze^{\xi }dx-%
\displaystyle\int_{\vartheta }A\nabla z\cdot \nabla \xi ze^{\xi }dx+%
\displaystyle\int_{\vartheta }gze^{\xi }dx+\displaystyle\frac{1}{2}%
\int_{\vartheta }\left\vert z\right\vert ^{2}\partial _{t}\xi e^{\xi }dx \\ 
& =-\displaystyle\int_{\vartheta }A\nabla z\cdot \nabla ze^{\xi }dx-%
\displaystyle\frac{1}{2}\int_{\vartheta }A\nabla \left( z^{2}\right) \cdot
\nabla \xi e^{\xi }dx+\displaystyle\int_{\vartheta }gze^{\xi }dx+%
\displaystyle\frac{1}{2}\int_{\vartheta }\left\vert z\right\vert
^{2}\partial _{t}\xi e^{\xi }dx \\ 
& =-\displaystyle\int_{\vartheta }A\nabla z\cdot \nabla ze^{\xi }dx+%
\displaystyle\frac{1}{2}\int_{\vartheta }\left\vert z\right\vert ^{2}\left(
\partial _{t}\xi +\nabla \cdot \left( A\nabla \xi \right) +A\nabla \xi \cdot
\nabla \xi \right) e^{\xi }dx+\displaystyle\int_{\vartheta }gze^{\xi }dx%
\text{ .}%
\end{array}%
\end{equation*}%
When $\varphi ={\partial }_{t}ze^{\xi }$, we have%
\begin{equation*}
\int_{\vartheta }\left\vert {\partial }_{t}z\right\vert ^{2}e^{\xi
}dx+\int_{\vartheta }A\nabla z\cdot {\partial }_{t}\nabla ze^{\xi
}dx+\int_{\vartheta }A\nabla z\cdot \nabla \xi {\partial }_{t}ze^{\xi
}dx=\int_{\vartheta }g{\partial }_{t}ze^{\xi }dx
\end{equation*}%
which implies%
\begin{equation*}
\begin{array}{ll}
& \quad \displaystyle\frac{1}{2}\frac{d}{dt}\int_{\vartheta }A\nabla z\cdot
\nabla ze^{\xi }dx-\displaystyle\frac{1}{2}\int_{\vartheta }\partial
_{t}A\nabla z\cdot \nabla ze^{\xi }dx \\ 
& =\displaystyle\int_{\vartheta }A\nabla z\cdot {\partial }_{t}\nabla
ze^{\xi }dx+\displaystyle\frac{1}{2}\int_{\vartheta }A\nabla z\cdot \nabla
z\partial _{t}\xi e^{\xi }dx \\ 
& =-\displaystyle\int_{\vartheta }\left\vert {\partial }_{t}z\right\vert
^{2}e^{\xi }dx-\displaystyle\int_{\vartheta }A\nabla z\cdot \nabla \xi {%
\partial }_{t}ze^{\xi }dx+\displaystyle\int_{\vartheta }g{\partial }%
_{t}ze^{\xi }dx+\displaystyle\frac{1}{2}\int_{\vartheta }A\nabla z\cdot
\nabla z\partial _{t}\xi e^{\xi }dx \\ 
& =-\displaystyle\int_{\vartheta }\left\vert {\partial }_{t}z+A\nabla z\cdot
\nabla \xi -\frac{1}{2}g\right\vert ^{2}e^{\xi }dx+\displaystyle\frac{1}{2}%
\int_{\vartheta }A\nabla z\cdot \nabla z\partial _{t}\xi e^{\xi }dx \\ 
& \quad +\displaystyle\int_{\vartheta }\left( {\partial }_{t}z-g\right)
A\nabla z\cdot \nabla \xi e^{\xi }dx+\displaystyle\int_{\vartheta
}\left\vert A\nabla z\cdot \nabla \xi \right\vert ^{2}e^{\xi }dx+%
\displaystyle\int_{\vartheta }\left\vert \frac{1}{2}g\right\vert ^{2}e^{\xi
}dx\text{ .}%
\end{array}%
\end{equation*}%
We compute $\displaystyle\int_{\vartheta }\left( {\partial }_{t}z-g\right)
A\nabla z\cdot \nabla \xi e^{\xi }dx$ by taking $\varphi =A\nabla z\cdot
\nabla \xi e^{\xi }$: One has with standard summation notations and $%
A=\left( A_{ij}\right) _{1\leq i,j\leq n}$ 
\begin{equation*}
\begin{array}{ll}
& \quad \displaystyle\int_{\vartheta }\left( {\partial }_{t}z-g\right)
A\nabla z\cdot \nabla \xi e^{\xi }dx+\displaystyle\int_{\vartheta
}\left\vert A\nabla z\cdot \nabla \xi \right\vert ^{2}e^{\xi }dx \\ 
& =-\displaystyle\int_{\vartheta }A\nabla z\cdot \nabla \left( A\nabla
z\cdot \nabla \xi e^{\xi }\right) dx+\displaystyle\int_{\partial \vartheta
}\left( A\nabla z\cdot \nu \right) \left( A\nabla z\cdot \nabla \xi \right)
e^{\xi }dx+\displaystyle\int_{\vartheta }\left\vert A\nabla z\cdot \nabla
\xi \right\vert ^{2}e^{\xi }dx \\ 
& =-\displaystyle\int_{\vartheta }A_{ij}\partial _{x_{j}}z\partial
_{x_{i}}A_{k\ell }\partial _{x_{\ell }}z\partial _{x_{k}}\xi e^{\xi }dx-%
\displaystyle\int_{\vartheta }A\nabla ^{2}\xi A\nabla z\cdot \nabla ze^{\xi
}dx \\ 
& \quad -\displaystyle\int_{\vartheta }A\nabla ^{2}zA\nabla z\cdot \nabla
\xi e^{\xi }dx+\displaystyle\int_{\partial \vartheta }\left( A\nabla z\cdot
\nu \right) \left( A\nabla z\cdot \nabla \xi \right) e^{\xi }dx\text{ .}%
\end{array}%
\end{equation*}%
But by one integration by parts%
\begin{equation*}
\begin{array}{ll}
& \quad -\displaystyle\int_{\vartheta }A\nabla ^{2}zA\nabla z\cdot \nabla
\xi e^{\xi }dx \\ 
& =-\displaystyle\frac{1}{2}\int_{\partial \vartheta }\left( A\nabla z\cdot
\nabla z\right) \left( A\nabla \xi \cdot \nu \right) e^{\xi }dx+\displaystyle%
\frac{1}{2}\int_{\vartheta }\partial _{x_{\ell }}A_{ij}\partial
_{x_{j}}zA_{k\ell }\partial _{x_{i}}z\partial _{x_{k}}\xi e^{\xi }dx \\ 
& \quad +\displaystyle\frac{1}{2}\int_{\vartheta }\left( A\nabla z\cdot
\nabla z\right) \nabla \cdot \left( A\nabla \xi \right) e^{\xi }dx+%
\displaystyle\frac{1}{2}\int_{\vartheta }\left( A\nabla z\cdot \nabla
z\right) \left( A\nabla \xi \cdot \xi \right) e^{\xi }dx\text{ .}%
\end{array}%
\end{equation*}%
The homogeneous Dirichlet boundary condition on $z$ implies $\nabla z=\nu
\partial _{\nu }z$ on $\partial \vartheta $. Therefore, one deduces%
\begin{equation*}
\begin{array}{ll}
\displaystyle\frac{1}{2}\frac{d}{dt}\int_{\vartheta }A\nabla z\cdot \nabla
ze^{\xi }dx & =\displaystyle\frac{1}{2}\int_{\vartheta }\partial _{t}A\nabla
z\cdot \nabla ze^{\xi }dx-\displaystyle\int_{\vartheta }\left\vert {\partial 
}_{t}z+A\nabla z\cdot \nabla \xi -\frac{1}{2}g\right\vert ^{2}e^{\xi }dx \\ 
& \quad -\displaystyle\int_{\vartheta }A\nabla ^{2}\xi A\nabla z\cdot \nabla
ze^{\xi }dx \\ 
& \quad +\displaystyle\frac{1}{2}\int_{\vartheta }\left( A\nabla z\cdot
\nabla z\right) \left( \partial _{t}\xi +\nabla \cdot \left( A\nabla \xi
\right) +A\nabla \xi \cdot \nabla \xi \right) e^{\xi }dx \\ 
& \quad -\displaystyle\int_{\vartheta }A_{ij}\partial _{x_{j}}z\partial
_{x_{i}}A_{k\ell }\partial _{x_{\ell }}z\partial _{x_{k}}\xi e^{\xi }dx+%
\displaystyle\frac{1}{2}\int_{\vartheta }\partial _{x_{\ell }}A_{ij}\partial
_{x_{j}}zA_{k\ell }\partial _{x_{i}}z\partial _{x_{k}}\xi e^{\xi }dx \\ 
& \quad +\displaystyle\frac{1}{2}\int_{\partial \vartheta }\left( A\nabla
z\cdot \nabla z\right) \left( A\nabla \xi \cdot \nu \right) e^{\xi }dx+%
\displaystyle\int_{\vartheta }\left\vert \frac{1}{2}g\right\vert ^{2}e^{\xi
}dx\text{ .}%
\end{array}%
\end{equation*}%
Now, we are able to compute $\displaystyle\frac{d}{dt}\frac{\int_{\vartheta
}A\nabla z\cdot \nabla ze^{\xi }dx}{\int_{\vartheta }\left\vert z\right\vert
^{2}e^{\xi }dx}$: One has%
\begin{equation*}
\begin{array}{ll}
& \quad \left( \displaystyle\int_{\vartheta }\left\vert z\right\vert
^{2}e^{\xi }dx\right) ^{2}\displaystyle\frac{d}{dt}\frac{\displaystyle%
\int_{\vartheta }A\nabla z\cdot \nabla ze^{\xi }dx}{\displaystyle%
\int_{\vartheta }\left\vert z\right\vert ^{2}e^{\xi }dx} \\ 
& =-\displaystyle2\int_{\vartheta }A\nabla ^{2}\xi A\nabla z\cdot \nabla
ze^{\xi }dx\displaystyle\int_{\vartheta }\left\vert z\right\vert ^{2}e^{\xi
}dx+\displaystyle\int_{\partial \vartheta }\left( A\nabla z\cdot \nabla
z\right) \left( A\nabla \xi \cdot \nu \right) e^{\xi }dx\displaystyle%
\int_{\vartheta }\left\vert z\right\vert ^{2}e^{\xi }dx \\ 
& \quad -\displaystyle2\int_{\vartheta }\left\vert {\partial }_{t}z+A\nabla
z\cdot \nabla \xi -\frac{1}{2}g\right\vert ^{2}e^{\xi }dx\displaystyle%
\int_{\vartheta }\left\vert z\right\vert ^{2}e^{\xi }dx \\ 
& \quad +\displaystyle2\int_{\vartheta }A\nabla z\cdot \nabla ze^{\xi
}dx\left( \displaystyle\int_{\vartheta }A\nabla z\cdot \nabla ze^{\xi }dx-%
\displaystyle\int_{\vartheta }gze^{\xi }dx\right) \\ 
& \quad +\displaystyle\int_{\vartheta }\partial _{t}A\nabla z\cdot \nabla
ze^{\xi }dx\displaystyle\int_{\vartheta }\left\vert z\right\vert ^{2}e^{\xi
}dx \\ 
& \quad +2\left( -\displaystyle\int_{\vartheta }A_{ij}\partial
_{x_{j}}z\partial _{x_{i}}A_{k\ell }\partial _{x_{\ell }}z\partial
_{x_{k}}\xi e^{\xi }dx+\displaystyle\frac{1}{2}\int_{\vartheta }\partial
_{x_{\ell }}A_{ij}\partial _{x_{j}}zA_{k\ell }\partial _{x_{i}}z\partial
_{x_{k}}\xi e^{\xi }dx\right) \displaystyle\int_{\vartheta }\left\vert
z\right\vert ^{2}e^{\xi }dx \\ 
& \quad +\displaystyle\int_{\vartheta }\left( A\nabla z\cdot \nabla z\right)
\left( \partial _{t}\xi +\nabla \cdot \left( A\nabla \xi \right) +A\nabla
\xi \cdot \nabla \xi \right) e^{\xi }dx\displaystyle\int_{\vartheta
}\left\vert z\right\vert ^{2}e^{\xi }dx \\ 
& \quad -\displaystyle\int_{\vartheta }A\nabla z\cdot \nabla ze^{\xi
}dx\left( \displaystyle\int_{\vartheta }\left\vert z\right\vert ^{2}\left(
\partial _{t}\xi +\nabla \cdot \left( A\nabla \xi \right) +A\nabla \xi \cdot
\nabla \xi \right) e^{\xi }dx\right) \\ 
& \quad +\displaystyle2\int_{\vartheta }\left\vert \frac{1}{2}g\right\vert
^{2}e^{\xi }dx\displaystyle\int_{\vartheta }\left\vert z\right\vert
^{2}e^{\xi }dx\text{ .}%
\end{array}%
\end{equation*}%
Notice that by Cauchy-Schwarz inequality, the contribution of the fourth and
fifth terms of the above becomes 
\begin{equation*}
\begin{array}{ll}
& \quad -\displaystyle\int_{\vartheta }\left\vert {\partial }_{t}z+A\nabla
z\cdot \nabla \xi -\frac{1}{2}g\right\vert ^{2}e^{\xi }dx\displaystyle%
\int_{\vartheta }\left\vert z\right\vert ^{2}e^{\xi }dx \\ 
& \quad +\displaystyle\int_{\vartheta }A\nabla z\cdot \nabla ze^{\xi
}dx\left( \displaystyle\int_{\vartheta }A\nabla z\cdot \nabla ze^{\xi }dx-%
\displaystyle\int_{\vartheta }gze^{\xi }dx\right) \\ 
& =-\displaystyle\int_{\vartheta }\left\vert {\partial }_{t}z+A\nabla z\cdot
\nabla \xi -\frac{1}{2}g\right\vert ^{2}e^{\xi }dx\displaystyle%
\int_{\vartheta }\left\vert z\right\vert ^{2}e^{\xi }dx \\ 
& \quad +\left( -\displaystyle\int_{\vartheta }\left( {\partial }%
_{t}z+A\nabla z\cdot \nabla \xi -\frac{1}{2}g\right) ze^{\xi }dx+%
\displaystyle\frac{1}{2}\int_{\vartheta }gze^{\xi }dx\right) \\ 
& \quad \times \left( -\displaystyle\int_{\vartheta }\left( {\partial }%
_{t}z+A\nabla z\cdot \nabla \xi -\frac{1}{2}g\right) ze^{\xi }dx-%
\displaystyle\frac{1}{2}\int_{\vartheta }gze^{\xi }dx\right) \\ 
& =-\displaystyle\int_{\vartheta }\left\vert {\partial }_{t}z+A\nabla z\cdot
\nabla \xi -\frac{1}{2}g\right\vert ^{2}e^{\xi }dx\displaystyle%
\int_{\vartheta }\left\vert z\right\vert ^{2}e^{\xi }dx \\ 
& \quad +\left( \displaystyle\int_{\vartheta }\left( {\partial }%
_{t}z+A\nabla z\cdot \nabla \xi -\frac{1}{2}g\right) ze^{\xi }dx\right)
^{2}-\left( \displaystyle\frac{1}{2}\int_{\vartheta }gze^{\xi }dx\right) ^{2}
\\ 
& \leq 0\text{ .}%
\end{array}%
\end{equation*}%
Therefore, one conclude that%
\begin{equation*}
\begin{array}{ll}
\displaystyle\frac{d}{dt}\frac{\displaystyle\int_{\vartheta }A\nabla z\cdot
\nabla ze^{\xi }dx}{\displaystyle\int_{\vartheta }\left\vert z\right\vert
^{2}e^{\xi }dx} & \leq \frac{-\displaystyle2\int_{\vartheta }A\nabla ^{2}\xi
A\nabla z\cdot \nabla ze^{\xi }dx}{\displaystyle\int_{\vartheta }\left\vert
z\right\vert ^{2}e^{\xi }dx}+\frac{\displaystyle\int_{\partial \vartheta
}\left( A\nabla z\cdot \nabla z\right) \left( A\nabla \xi \cdot \nu \right)
e^{\xi }dx}{\displaystyle\int_{\vartheta }\left\vert z\right\vert ^{2}e^{\xi
}dx} \\ 
& \quad +\frac{\displaystyle\int_{\vartheta }\left\vert g\right\vert
^{2}e^{\xi }dx}{\displaystyle\int_{\vartheta }\left\vert z\right\vert
^{2}e^{\xi }dx}+\frac{\displaystyle\int_{\vartheta }\partial _{t}A\nabla
z\cdot \nabla ze^{\xi }dx}{\displaystyle\int_{\vartheta }\left\vert
z\right\vert ^{2}e^{\xi }dx} \\ 
& \quad +\frac{-\displaystyle2\int_{\vartheta }A_{ij}\partial
_{x_{j}}z\partial _{x_{i}}A_{k\ell }\partial _{x_{\ell }}z\partial
_{x_{k}}\xi e^{\xi }dx+\displaystyle\int_{\vartheta }\partial _{x_{\ell
}}A_{ij}\partial _{x_{j}}zA_{k\ell }\partial _{x_{i}}z\partial _{x_{k}}\xi
e^{\xi }dx}{\displaystyle\int_{\vartheta }\left\vert z\right\vert ^{2}e^{\xi
}dx} \\ 
& \quad +\frac{\displaystyle\int_{\vartheta }A\nabla z\cdot \nabla z\left(
\partial _{t}\xi +\nabla \cdot \left( A\nabla \xi \right) +A\nabla \xi \cdot
\nabla \xi \right) e^{\xi }dx}{\displaystyle\int_{\vartheta }\left\vert
z\right\vert ^{2}e^{\xi }dx} \\ 
& \quad -\frac{\displaystyle\int_{\vartheta }A\nabla z\cdot \nabla ze^{\xi
}dx}{\displaystyle\int_{\vartheta }\left\vert z\right\vert ^{2}e^{\xi }dx}%
\times \frac{\displaystyle\int_{\vartheta }\left\vert z\right\vert
^{2}\left( \partial _{t}\xi +\nabla \cdot \left( A\nabla \xi \right)
+A\nabla \xi \cdot \nabla \xi \right) e^{\xi }dx}{\displaystyle%
\int_{\vartheta }\left\vert z\right\vert ^{2}e^{\xi }dx}\text{ .}%
\end{array}%
\end{equation*}

\bigskip

\bigskip

\subsection{Proof of Lemma \protect\ref{lemma4.7}}

\bigskip

Let $0<\varepsilon <\rho /2$ and $\phi \in C_{0}^{\infty }\left( B_{\rho
}\right) $ be such that $0\leq \phi \leq 1$, $\phi =1$ on $\left\{ x;d\left(
x,x_{0}\right) \leq \rho -\varepsilon \right\} $. We multiply the equation $%
\partial _{t}u-\nabla \cdot \left( A\nabla u\right) =0$ by $\phi
^{2}ue^{-d\left( x,x_{0}\right) ^{2}/h}$ where $h>0$ and integrate over $%
\Omega \cap B_{\rho }$. We get by one integration by parts 
\begin{equation*}
\frac{1}{2}\frac{d}{dt}\int_{\Omega \cap B_{\rho }}\left\vert \phi
u\right\vert ^{2}e^{-d\left( x,x_{0}\right) ^{2}/h}dx+\int_{\Omega \cap
B_{\rho }}A\nabla u\cdot \nabla \left( \phi ^{2}ue^{-d\left( x,x_{0}\right)
^{2}/h}\right) dx=0\text{ .}
\end{equation*}%
But, $A\nabla u\cdot \nabla \left( \phi ^{2}ue^{-d^{2}/h}\right) =\left[
2\phi uA\nabla \phi \cdot \nabla u+\phi ^{2}A\nabla u\cdot \nabla u+\phi
^{2}u\left( -\frac{2d\nabla d}{h}\right) \cdot A\nabla u\right] e^{-d^{2}/h}$%
. Therefore, by Cauchy-Schwarz inequality, it follows that 
\begin{equation*}
\begin{array}{ll}
& \quad \displaystyle\frac{1}{2}\frac{d}{dt}\int_{\Omega \cap B_{\rho
}}\left\vert \phi u\right\vert ^{2}e^{-d\left( x,x_{0}\right) ^{2}/h}dx+%
\displaystyle\int_{\Omega \cap B_{\rho }}\phi ^{2}A\nabla u\cdot \nabla
ue^{-d\left( x,x_{0}\right) ^{2}/h}dx \\ 
& \leq \displaystyle\int_{\Omega \cap B_{\rho }}\phi ^{2}A\nabla u\cdot
\nabla ue^{-d\left( x,x_{0}\right) ^{2}/h}dx+\displaystyle\frac{1}{2}%
\int_{\Omega \cap B_{\rho }}4A\nabla \phi \cdot \nabla \phi \left\vert
u\right\vert ^{2}e^{-d\left( x,x_{0}\right) ^{2}/h}dx \\ 
& \quad +\displaystyle\frac{1}{2}\int_{\Omega \cap B_{\rho }}\frac{4d^{2}}{%
h^{2}}A\nabla d\cdot \nabla d\left\vert \phi u\right\vert ^{2}e^{-d\left(
x,x_{0}\right) ^{2}/h}dx\text{ .}%
\end{array}%
\end{equation*}%
Thus, with the fact that $A_{T}\left( x\right) \nabla d\left( x,x_{0}\right)
\cdot \nabla d\left( x,x_{0}\right) =1$, one get for some constant $C_{A}>0$ 
\begin{equation*}
\begin{array}{ll}
& \quad \displaystyle\frac{d}{dt}\int_{\Omega \cap B_{\rho }}\left\vert \phi
u\right\vert ^{2}e^{-d\left( x,x_{0}\right) ^{2}/h}dx-\displaystyle\frac{%
\rho ^{2}}{h^{2}}C_{A}\int_{\Omega \cap B_{\rho }}\left\vert \phi
u\right\vert ^{2}e^{-d\left( x,x_{0}\right) ^{2}/h}dx \\ 
& \leq C_{A}e^{-\frac{\left( \rho -\varepsilon \right) ^{2}}{h}}\displaystyle%
\int_{\Omega \cap B_{\rho }}\left\vert u\left( x,t\right) \right\vert ^{2}dx%
\text{ .}%
\end{array}%
\end{equation*}%
Then we have, 
\begin{equation*}
\begin{array}{ll}
\displaystyle\int_{\Omega \cap B_{\rho }}\left\vert \phi u\left( \cdot
,T\right) \right\vert ^{2}e^{-d\left( x,x_{0}\right) ^{2}/h}dx & \leq e^{%
\frac{\rho ^{2}}{h^{2}}C_{A}\left( T-t\right) }\displaystyle\int_{\Omega
\cap B_{\rho }}\left\vert \phi u\left( \cdot ,t\right) \right\vert
^{2}e^{-d\left( x,x_{0}\right) ^{2}/h}dx \\ 
& \quad +C_{A}e^{\frac{\rho ^{2}}{h^{2}}C_{A}\left( T-t\right) }e^{-\frac{%
\left( \rho -\varepsilon \right) ^{2}}{h}}\displaystyle\int_{t}^{T}%
\displaystyle\int_{\Omega \cap B_{\rho }}\left\vert u\right\vert ^{2}dxds%
\end{array}%
\end{equation*}%
which gives%
\begin{equation*}
\begin{array}{ll}
\displaystyle\int_{\Omega \cap B_{\rho -2\varepsilon }}\left\vert u\left(
x,T\right) \right\vert ^{2}dx & \leq e^{\frac{\rho ^{2}}{h^{2}}C_{A}\left(
T-t\right) }e^{\frac{\left( \rho -2\varepsilon \right) ^{2}}{h}}\displaystyle%
\int_{\Omega \cap B_{\rho }}\left\vert u\left( x,t\right) \right\vert ^{2}dx
\\ 
& \quad +C_{A}e^{\frac{\rho ^{2}}{h^{2}}C_{A}\left( T-t\right) }e^{-\frac{%
\left( \rho -\varepsilon \right) ^{2}}{h}}e^{\frac{\left( \rho -2\varepsilon
\right) ^{2}}{h}}\displaystyle\int_{t}^{T}\displaystyle\int_{\Omega \cap
B_{\rho }}\left\vert u\right\vert ^{2}dxds\text{ .}%
\end{array}%
\end{equation*}%
Let $T/2<T-\delta h\leq t\leq T$, it yields 
\begin{equation*}
\begin{array}{ll}
\displaystyle\int_{\Omega \cap B_{\rho -2\varepsilon }}\left\vert u\left(
x,T\right) \right\vert ^{2}dx & \leq e^{CT}e^{\frac{\rho ^{2}}{h}\delta
C_{A}}e^{\frac{\left( \rho -2\varepsilon \right) ^{2}}{h}}\displaystyle%
\int_{\Omega \cap B_{\rho }}\left\vert u\left( x,t\right) \right\vert ^{2}dx
\\ 
& \quad +C_{A}e^{\frac{\rho ^{2}}{h}\delta C_{A}}e^{-\frac{\left( \rho
-\varepsilon \right) ^{2}}{h}}e^{\frac{\left( \rho -2\varepsilon \right) ^{2}%
}{h}}\displaystyle\int_{T-\delta h}^{T}\displaystyle\int_{\Omega \cap
B_{\rho }}\left\vert u\right\vert ^{2}dxds\text{ .}%
\end{array}%
\end{equation*}%
Choose 
\begin{equation*}
\delta =\frac{1}{C_{A}}\frac{\varepsilon \left( 2\rho -3\varepsilon \right) 
}{2\rho ^{2}}
\end{equation*}%
that is $\delta C_{A}=\frac{1}{2}\frac{\left( \rho -\varepsilon \right)
^{2}-\left( \rho -2\varepsilon \right) ^{2}}{\rho ^{2}}\in \left( 0,1/8%
\right] $ in order that $\rho ^{2}\delta C_{A}-\left( \rho -\varepsilon
\right) ^{2}+\left( \rho -2\varepsilon \right) ^{2}<0$. Therefore, we get%
\begin{equation*}
\begin{array}{ll}
\displaystyle\int_{\Omega \cap B_{\rho -2\varepsilon }}\left\vert u\left(
x,T\right) \right\vert ^{2}dx & \leq e^{\frac{\left( \rho -\varepsilon
\right) ^{2}+\left( \rho -2\varepsilon \right) ^{2}}{2h}}\displaystyle%
\int_{\Omega \cap B_{\rho }}\left\vert u\left( x,t\right) \right\vert ^{2}dx
\\ 
& \quad +C_{A}e^{\frac{-\left( \rho -\varepsilon \right) ^{2}+\left( \rho
-2\varepsilon \right) ^{2}}{2h}}\displaystyle\int_{T-\delta h}^{T}%
\displaystyle\int_{\Omega \cap B_{\rho }}\left\vert u\right\vert ^{2}dxdt \\ 
& \leq e^{\frac{\left( \rho -\varepsilon \right) ^{2}+\left( \rho
-2\varepsilon \right) ^{2}}{2h}}\displaystyle\int_{\Omega \cap B_{\rho
}}\left\vert u\left( x,t\right) \right\vert ^{2}dx \\ 
& \quad +C_{A}e^{\frac{-\left( \rho -\varepsilon \right) ^{2}+\left( \rho
-2\varepsilon \right) ^{2}}{2h}}\displaystyle\int_{\Omega }\left\vert
u\left( x,0\right) \right\vert ^{2}dx%
\end{array}%
\end{equation*}%
where in the last line we used $\delta h<$max$\left( 1,T/2\right) $. Now,
choose $h$ such that both $\delta h<$max$\left( 1,T/2\right) $ and 
\begin{equation*}
\left( 1+C_{A}\right) e^{\frac{-\left( \rho -\varepsilon \right) ^{2}+\left(
\rho -2\varepsilon \right) ^{2}}{2h}}\int_{\Omega }\left\vert u\left(
x,0\right) \right\vert ^{2}dx\leq \dfrac{1}{e}\int_{\Omega \cap B_{\rho
-2\varepsilon }}\left\vert u\left( x,T\right) \right\vert ^{2}dx\text{ .}
\end{equation*}%
With such choice, one has%
\begin{equation*}
\left( 1-\frac{1}{e}\right) \int_{\Omega \cap B_{\rho -2\varepsilon
}}\left\vert u\left( x,T\right) \right\vert ^{2}dx\leq e^{\frac{\left( \rho
-\varepsilon \right) ^{2}+\left( \rho -2\varepsilon \right) ^{2}}{2h}%
}\int_{\Omega \cap B_{\rho }}\left\vert u\left( x,t\right) \right\vert ^{2}dx
\end{equation*}%
and moreover,%
\begin{equation*}
\int_{\Omega }\left\vert u\left( x,0\right) \right\vert ^{2}dx\leq e^{\frac{%
\left( \rho -\varepsilon \right) ^{2}}{h}}\int_{\Omega \cap B_{\rho
}}\left\vert u\left( x,t\right) \right\vert ^{2}dx
\end{equation*}%
for any $T/2<T-\delta h\leq t\leq T$. Such $h$ exists by choosing 
\begin{equation*}
h=\frac{\varepsilon \left( 2\rho -3\varepsilon \right) /2}{\text{ln}\left( K%
\frac{\left( 1+C_{A}\right) \displaystyle\int_{\Omega }\left\vert u\left(
x,0\right) \right\vert ^{2}dx}{\dfrac{1}{e}\displaystyle\int_{\Omega \cap
B_{\rho -2\varepsilon }}\left\vert u\left( x,T\right) \right\vert ^{2}dx}%
\right) }\text{ with }K=e^{\varepsilon \frac{\left( 2\rho -3\varepsilon
\right) }{2}\left( \frac{2}{T}+1\right) \delta }\text{ .}
\end{equation*}%
Clearly, $\delta h<T/2$ and $\delta h\leq 1$. We conclude that for any $%
T/2\leq T-\theta \leq t\leq T$ 
\begin{equation*}
\frac{\displaystyle\int_{\Omega }\left\vert u\left( x,0\right) \right\vert
^{2}dx}{\displaystyle\int_{\Omega \cap B_{\rho }}\left\vert u\left(
x,t\right) \right\vert ^{2}dx}\leq e^{\frac{1}{C_{A}}\frac{\varepsilon
\left( 2\rho -3\varepsilon \right) \left( \rho -\varepsilon \right) ^{2}}{%
2\rho ^{2}}\frac{1}{\theta }}
\end{equation*}%
with 
\begin{equation*}
\frac{1}{\theta }=C_{A}\frac{4\rho ^{2}}{\varepsilon ^{2}\left( 2\rho
-3\varepsilon \right) ^{2}}\text{ln}\left( e\left( 1+C_{A}\right) e^{\left( 
\frac{2}{T}+1\right) \frac{1}{C_{A}}\frac{\varepsilon ^{2}\left( 2\rho
-3\varepsilon \right) ^{2}}{4\rho ^{2}}}\frac{\displaystyle\int_{\Omega
}\left\vert u\left( x,0\right) \right\vert ^{2}dx}{\displaystyle\int_{\Omega
\cap B_{\rho -2\varepsilon }}\left\vert u\left( x,T\right) \right\vert ^{2}dx%
}\right) \text{ .}
\end{equation*}%
This completes the proof.

\bigskip

Remark .- When $A$ is time-independent, then $C_{A}=4$max$\left(
1,\left\Vert A\nabla \phi \cdot \nabla \phi \right\Vert _{L^{\infty }\left(
\Omega \right) }\right) $.

\bigskip

\bigskip

\section*{Appendix}

\bigskip

This appendix is devoted to the proof of Proposition \ref{proposition3.2}
and of inequality (\ref{2.2}).

\bigskip

\subsection*{Trace estimate for $f$ (proof of Proposition \protect\ref%
{proposition3.2})}

\bigskip

Denote $\left( \partial \Omega \times \mathbb{S}^{\text{d}-1}\right)
_{+}=\left\{ \left( x,v\right) \in \partial \Omega \times \mathbb{S}^{\text{d%
}-1};v\cdot \vec{n}_{x}\geq 0\right\} $. First, multiplying both sides of
the first line of (\ref{1.1}) by $\eta f\left\vert f\right\vert ^{\eta -2}$
and integrating over $\Omega \times \mathbb{S}^{\text{d}-1}\times \left(
0,T\right) $, one has the following a priori estimate for any $\eta \geq 2$%
\begin{equation*}
\int_{0}^{T}\int_{\left( \partial \Omega \times \mathbb{S}^{\text{d}%
-1}\right) _{+}}v\cdot \vec{n}_{x}\left\vert f\right\vert ^{\eta }dxdvdt\leq
\epsilon \frac{2}{\eta }\int_{\Omega \times \mathbb{S}^{\text{d}%
-1}}\left\vert f_{0}\right\vert ^{\eta }dxdv\text{ .}
\end{equation*}%
Secondly, one uses H\"{o}lder inequality to get 
\begin{equation*}
\begin{array}{ll}
& \quad \displaystyle\int_{0}^{T}\int_{\partial \Omega \times \mathbb{S}^{%
\text{d}-1}}\left\vert f\right\vert ^{2}dxdvdt \\ 
& \leq \left( \displaystyle\int_{0}^{T}\int_{\left( \partial \Omega \times 
\mathbb{S}^{\text{d}-1}\right) _{+}}\dfrac{dxdvdt}{\left( v\cdot \vec{n}%
_{x}\right) ^{\frac{1}{p-1}}}\right) ^{\frac{p-1}{p}}\left( \displaystyle%
\int_{0}^{T}\int_{\left( \partial \Omega \times \mathbb{S}^{\text{d}%
-1}\right) _{+}}v\cdot \vec{n}_{x}\left\vert f\right\vert ^{2p}dxdvdt\right)
^{\frac{1}{p}}\text{ .}%
\end{array}%
\end{equation*}%
But%
\begin{equation*}
\int_{\left( \partial \Omega \times \mathbb{S}^{\text{d}-1}\right) _{+}}%
\dfrac{dxdv}{\left( v\cdot \vec{n}_{x}\right) ^{\frac{1}{p-1}}}\leq C\frac{%
p-1}{p-2}\text{ for any }p>2\text{ .}
\end{equation*}%
Hence, as soon as $p>2$, one get the desired estimate 
\begin{equation*}
\left\Vert f\right\Vert _{L^{2}\left( \partial \Omega \times \mathbb{S}^{%
\text{d}-1}\times \left( 0,T\right) \right) }\leq CT^{\frac{p-1}{2p}%
}\epsilon ^{\frac{1}{2p}}C_{p}\Vert f_{0}\Vert _{L^{2p}(\Omega \times 
\mathbb{S}^{\text{d}-1})}
\end{equation*}%
where $C_{p}=\left( \frac{p-1}{p-2}\right) ^{\frac{p-1}{2p}}\left( \frac{1}{p%
}\right) ^{\frac{1}{2p}}$ and $C>0$ only depends on $\left( \Omega ,\text{d}%
\right) $.

\bigskip

\subsection*{Backward estimate for diffusion equations (proof of (\protect
\ref{2.2}))}

\bigskip

Classical energy identities for our parabolic equation are:%
\begin{equation*}
\frac{1}{2}\frac{d}{dt}\int_{\Omega}\left\vert u\right\vert
^{2}dx+\int_{\Omega}\frac{1}{\text{d}a}\left\vert \nabla u\right\vert
^{2}dx=0\text{ ,}
\end{equation*}%
\begin{equation*}
\frac{1}{2}\frac{d}{dt}\int_{\Omega}\frac{1}{\text{d}a}\left\vert
\nabla\varphi\right\vert ^{2}dx+\int_{\Omega}\left\vert u\right\vert ^{2}dx=0%
\text{ ,}
\end{equation*}
where $\varphi\left( \cdot,t\right) \in H_{0}^{1}\left( \Omega\right) $
solves $-\nabla\cdot\left( \frac{1}{\text{d}a}\nabla\varphi\left(
\cdot,t\right) \right) =u\left( \cdot,t\right) $ in $\Omega$. Now, one can
easily check with $y\left( t\right) =\displaystyle\int_{\Omega}\frac {1}{%
\text{d}a\left( x\right) }\left\vert \nabla\varphi\left( x,t\right)
\right\vert ^{2}dx$ and $N\left( t\right) =\frac{\displaystyle\int_{\Omega
}\left\vert u\left( x,t\right) \right\vert ^{2}dx}{\displaystyle\int
_{\Omega}\frac{1}{\text{d}a\left( x\right) }\left\vert \nabla\varphi\left(
x,t\right) \right\vert ^{2}dx}$ that 
\begin{equation*}
\left\{ 
\begin{array}{ll}
\frac{1}{2}y^{\prime}\left( t\right) +N\left( t\right) y\left( t\right) =0 & 
\\ 
N^{\prime}\left( t\right) \leq0\text{ .} & 
\end{array}
\right.
\end{equation*}
By solving such differential inequalities, one obtain 
\begin{equation*}
\int_{\Omega}\frac{1}{\text{d}a\left( x\right) }\left\vert \nabla
\varphi\left( x,0\right) \right\vert ^{2}dx\leq e^{2TN\left( 0\right)
}\int_{\Omega}\frac{1}{\text{d}a\left( x\right) }\left\vert \nabla
\varphi\left( x,T\right) \right\vert ^{2}dx
\end{equation*}
which implies 
\begin{equation*}
\left\Vert u\left( \cdot,T\right) \right\Vert _{H^{-1}\left( \Omega\right)
}^{2}\leq\frac{c_{max}}{c_{min}}e^{2T\frac{\left\Vert u\left( \cdot,0\right)
\right\Vert _{L^{2}\left( \Omega\right) }^{2}}{\text{d}c_{min}\left\Vert
u\left( \cdot,0\right) \right\Vert _{H^{-1}\left( \Omega\right) }^{2}}%
}\left\Vert u\left( \cdot,T\right) \right\Vert _{H^{-1}\left( \Omega\right)
}^{2}\text{ .}
\end{equation*}
One conclude that 
\begin{equation*}
\left\Vert u\left( \cdot,0\right) \right\Vert _{L^{2}\left( \Omega\right)
}\leq ce^{cT\frac{\left\Vert u\left( \cdot,0\right) \right\Vert
_{L^{2}\left( \Omega\right) }^{2}}{\left\Vert u\left( \cdot,0\right)
\right\Vert _{H^{-1}\left( \Omega\right) }^{2}}}\left\Vert u\left(
\cdot,T\right) \right\Vert _{L^{2}\left( \Omega\right) }\text{ .}
\end{equation*}

\bigskip

\bigskip

\bigskip

\bigskip

\bigskip

\bigskip


\bigskip

\bigskip

\bigskip

\begin{thebibliography}{AEWZ}
\bibitem[AEWZ]{AEWZ} J. Apraiz, L. Escauriaza, G. Wang and C. Zhang,
Observability inequalities and measurable sets. Journal of the European
Mathematical Society 16 (2014), no. 11, 2433--2475.

\bibitem[BR]{BR} G. Bal and L. Ryzhik, Diffusion approximation of radiative
transfer problems with interfaces. SIAM J. Appl. Math. 60 (2000), no. 6,
1887--1912.

\bibitem[B]{B} C. Bardos, Probl\`{e}mes aux limites pour les \'{e}quations
aux d\'{e}riv\'{e}es partielles du premier ordre \`{a} coefficients r\'{e}%
els; th\'{e}or\`{e}mes d'approximation; applications \`{a} l'\'{e}quations
de tranport. Ann. Sci. Ecole Norm. Sup. 3 (1970), 185--233.

\bibitem[BBGS]{BBGS} C. Bardos, E. Bernard, F. Golse and R. Sentis, The
diffusion approximation for the linear Boltzmann equation with vanishing
scattering coefficient, Communications in Mathematical Sciences 13 (2015),
no. 3, 641--671.

\bibitem[BGPS]{BGPS} C. Bardos, F. Golse, B. Perthame and R. Sentis, The
nonaccretive radiative transfer equations: existence of solutions and
Rosseland approximation. J. Funct. Anal. 77 (1988), no. 2, 434--460.

\bibitem[BSS]{BSS} C. Bardos, R. Santos and R. Sentis, Diffusion
approximation and computation of the critical size. Trans. of the AMS. 284
(1984), 617--649.

\bibitem[BT]{BT} C. Bardos and L. Tartar, Sur l'unicit\'{e} retrograde des 
\'{e}quations paraboliques et quelques questions voisines. Arch. Rational
Mech. Anal. 50 (1973), 10--25.

\bibitem[DL]{DL} R. Dautray and J.-L. Lions, Analyse Math\'{e}matique et
Calcul Num\'{e}rique pour les Sciences et les Techniques, Masson, 1985.

\bibitem[DJP]{DJP} D. Del Santo, C. J\"{a}h and M. Paicu, Backward
uniqueness for parabolic operators with non-Lipschitz coefficients. Osaka J.
Math. 52 (2015), no. 3, 793--815.

\bibitem[EFV]{EFV} L. Escauriaza, F.J. Fernandez and S. Vessella, Doubling
properties of caloric functions, Appl. Anal. 85 (2006), 205--223.

\bibitem[FG]{FG} E. Fern\'{a}ndez-Cara and S. Guerrero, Global Carleman
inequalities for parabolic systems and application to controllability. SIAM
J. Control Optim. 45 (2006), no. 4, 1395--1446.

\bibitem[FZ]{FZ} E. Fern\'{a}ndez-Cara and E. Zuazua, Null and approximate
controllability for weakly blowing up semilinear heat equations. Annales de
l'Institut Henri Poincar\'{e}, Analyse Non Lin\'{e}aire 17 (2000), 583--616.

\bibitem[FI]{FI} A.V. Fursikov and O.Yu. Imanuvilov, Controllability of
evolution equations. Lecture Notes Series, 34. Seoul National University,
Research Institute of Mathematics, Global Analysis Research Center, Seoul,
1996.

\bibitem[IM]{IM} S. Ibrahim and M. Majdoub, Solutions globales de l'\'{e}%
quation des ondes semi-lin\'{e}aire critique \`{a} coefficients variables.
Bull. Soc. Math. France 131 (2003), no. 1, 1--22.

\bibitem[I]{I} V. Isakov, Inverse Problems for Partial Differential
Equations, Second Edition, Springer, New York, 2006.

\bibitem[LK]{LK} E.W. Larsen and J.B. Keller, Asymptotic solution of neutron
transport problems for small mean free paths. J. Math. Phys. 15 (1974),
75--81.

\bibitem[LR]{LR} G. Lebeau and L. Robbiano, Contr\^{o}le exact de l'\'{e}%
quation de la chaleur. Comm. Partial Differential Equations 20 (1995), no.
1-2, 335--356.

\bibitem[LRL]{LRL} J. Le Rousseau and G. Lebeau, On Carleman estimates for
elliptic and parabolic operators. Applications to unique continuation and
control of parabolic equations. ESAIM: Control, Optimisation and Calculus of
Variations 18 (2012), 712--747.

\bibitem[LeRR]{LeRR} J. Le Rousseau and L. Robbiano, Carleman estimate for
elliptic operators with coefficients with jumps at an interface in arbitrary
dimension and application to the null controllability of linear parabolic
equations. Arch. Ration. Mech. Anal. 195 (2010), no. 3, 953--990.

\bibitem[LRR]{LRR} J. Le Rousseau and L. Robbiano, Local and global Carleman
estimates for parabolic operators with coefficients with jumps at
interfaces. Invent. Math. 183 (2011), no. 2, 245--336.

\bibitem[LM]{LM} J.-L. Lions and B. Malgrange, Sur l'unicit\'{e} r\'{e}%
trograde dans les probl\`{e}mes mixtes paraboliques. Math. Scand. 8 (1960),
277--286.

\bibitem[P]{P} L. Payne, Improperly Posed Problems in Partial Differential
Equations, Regional Conference Series in Applied Mathematics, Vol. 22 ,
SIAM, 1975.

\bibitem[PWa]{PWa} K.D. Phung and G. Wang, Quantitative unique continuation
for the semilinear heat equation in a convex domain. J. Funct. Anal. 259
(2010), no. 5, 1230--1247.

\bibitem[PW]{PW} K.D. Phung and G. Wang, An observability estimate for
parabolic equations from a measurable set in time and its applications.
Journal of the European Mathematical Society 15 (2013), no. 2, 681--703.

\bibitem[V]{V} S. Vessella, Unique continuation properties and quantitative
estimates of unique continuation for parabolic equations. Handbook of
differential equations: evolutionary equations. Vol. 5, 421-500, Handb.
Differ. Equ., Elsevier/North-Holland, Amsterdam, 2009.
\end{thebibliography}
\end{document}